\documentclass[11pt]{amsart}
\usepackage{geometry}                
\usepackage{url}
\usepackage{graphicx}
\usepackage{amssymb,mathtools}
\usepackage{amsthm}
\usepackage{epstopdf}
\usepackage{listings}
\usepackage{color}
\usepackage{bbold}
\theoremstyle{definition}

\DeclareMathOperator{\trace}{trace}
\title{Computational Bifurcation Analysis}
\author{Harry Dankowicz \& Jan Sieber}

\newcommand{\hiddenB}{white}
\newcommand{\hiddenA}{green}

\definecolor{grey1}{rgb}{0.0,0.0,0.0}
\definecolor{grey2}{rgb}{0.2,0.2,0.2}
\definecolor{grey3}{rgb}{0.4,0.4,0.4}
\definecolor{grey4}{rgb}{0.5,0.5,0.5}

\lstdefinestyle{highlight-fonts}{
	basicstyle=\ttfamily\mdseries\footnotesize,
	commentstyle=\ttfamily\mdseries\itshape\footnotesize,
	moredelim=[l][\color{\hiddenA}\ttfamily\mdseries\itshape\footnotesize]{\%\#},
	moredelim=[l][\color{\hiddenB}\ttfamily\mdseries\itshape\footnotesize]{\%!},
	keywordstyle=\ttfamily\mdseries\footnotesize,
	keywordstyle={[2]\ttfamily\bfseries\footnotesize},
}

\lstdefinestyle{highlight-colors}{
	basicstyle=\color{grey1}\ttfamily\mdseries\footnotesize,
	commentstyle=\color{grey2}\ttfamily\mdseries\itshape\footnotesize,
	moredelim=[l][\color{\hiddenA}\ttfamily\mdseries\itshape\footnotesize]{\%\#},
	moredelim=[l][\color{\hiddenB}\ttfamily\mdseries\itshape\footnotesize]{\%!},
	keywordstyle=\color{grey2}\ttfamily\mdseries\footnotesize,
	keywordstyle={[2]\color{grey2}\ttfamily\bfseries\footnotesize},
}

\lstdefinestyle{gobble1}{
	xleftmargin=0.3em,
}

\lstdefinestyle{gobble2}{
	xleftmargin=-0.9em,
}

\lstdefinestyle{gobble3}{
	xleftmargin=-2.1em,
}

\lstdefinelanguage{coco}[]{Matlab}{
	keywords={},
	morekeywords={},
	morekeywords=[2]{classdef,properties,private,protected,public,%
Access,Static,methods,if,function,end,for,while,else,elseif,switch,%
case,otherwise,do,repeat,until},
	otherkeywords={{,end},:end,(end,end),\{end,end\},[end,end],/private,private/},
	morecomment=[l]{},
	basicstyle=\ttfamily\mdseries\footnotesize,
	commentstyle=\ttfamily\mdseries\footnotesize,
	moredelim=[l][\color{\hiddenA}\ttfamily\mdseries\footnotesize]{\%\#},
	moredelim=[l][\color{\hiddenB}\ttfamily\mdseries\footnotesize]{\%!},
	keywordstyle=\ttfamily\mdseries\footnotesize,
	keywordstyle={[2]\ttfamily\mdseries\footnotesize},
	numbers=none,
	numberstyle=\scriptsize,
	numbersep=1em,
	breaklines=false,
	breakatwhitespace=true,
	breakindent=2.5em,
	showlines=false,
	lineskip=-0.2ex,
	frame=none,
	fontadjust=true,
	columns=[c]fixed,
	basewidth={0.575em,0.45em},
	fontadjust=true,
	keepspaces=true,
	tabsize=3,
	showstringspaces=false,
	aboveskip=1.5\medskipamount,
	belowskip=1.5\medskipamount,
	xleftmargin=0.925em,
	xrightmargin=0em,
	rangeprefix={\%!},
	includerangemarker=false,
	belowcaptionskip=\bigskipamount
}
\lstdefinelanguage{coco-highlight-fonts}[]{coco}{
	style=highlight-fonts,
}

\lstdefinelanguage{coco-highlight-colors}[]{coco}{
	style=highlight-colors,
}

\lstdefinelanguage{coco-highlight}[]{coco-highlight-colors}{}

\newcommand{\mcode}[1]{{\lstinline[language=coco-highlight,basewidth={0.6em,0.45em}]|#1|}}

\newcommand{\tran}{\mathsf{T}}
\renewcommand{\i}{\mathrm{i}}
\newcommand{\e}{\mathrm{e}}

\newcommand{\jsq}[1]{\textbf{\color{red}[JS: #1]}}
\begin{document}
\maketitle
\begin{flushright}
    \emph{Dedicated to the memory of Claudia Wulff}
\end{flushright}
\section{Introduction}
Bifurcation analysis collects techniques for characterizing the dependence of certain classes of solutions of a dynamical system on variations in problem parameters. Common solution classes of interest include equilibria and periodic orbits, the number and stability of which may vary as parameters vary. Continuation techniques generate continuous families of such solutions in the combined state and parameter space, e.g., curves (branches) of periodic orbits or surfaces of equilibria. Their advantage over simulation-based approaches is the ability to map out such families independently of the dynamic stability of the equilibria or periodic orbits. Bifurcation diagrams represent families of equilibria and periodic orbits as curves or surfaces in appropriate coordinate systems. Special points, such as bifurcations, are often highlighted in such diagrams.

Proficiency with bifurcation analysis requires a high degree of competence with both theory and computation. Optimal mileage is obtained through a judicious interplay between the two, seeing each as an extension of the other. Where one excels, the other may step back momentarily, only to come roaring back when the time is right. The analyst wields each of a multitude of theoretical and algorithmic approaches as so many arrows in a quiver, judging by the direction and strength of the wind, and by the distance to the target, the right tool for each moment.

With this perspective in mind, complete automation of analysis is neither desirable nor achievable. Undesirable, since it removes the analyst from the need to align technique against objective. Unachievable, since it is in the nature of the beast that the problems that are worthy of particular struggle are also often those that do not conform easily with established knowledge. Nevertheless, many tasks can be partially automated and it is certainly in one's best interest to avail oneself of such partial automation when it is wielded as an element of a systematic strategy and informed by an understanding of what can and cannot be accomplished.

This article provides an illustration of this paradigm of synergy between theoretical derivations and computational analysis for several characteristic examples of bifurcation analysis in commonly encountered classes of problems. General theoretical principles are deduced from these illustrations and collected for the reader's subsequent reference.

A degree of the promised partial automation is provided in this article by the \textsc{coco} package of Matlab-compatible software algorithms. As a companion to the content included in the chapter, a code repository on GitHub contains executable scripts to reproduce all the numerical results reported below.

\section{Illustration: Continuous Stirred Tank Reactor}
\label{sec:cstr}
As a first illustration of a general methodology for computational bifurcation analysis, the reader is invited to consider the planar dynamical system described by the two coupled ordinary differential equations
\begin{equation}
\label{eq:cstr:model}
\dot{x}=f_x(x,y):=(1-x)e^\frac{y}{1+\beta y}-\frac{x}{\delta},\,\gamma\dot{y}=\gamma f_y(x,y):=(1-x)e^\frac{y}{1+\beta y}-\frac{y}{\sigma}
\end{equation}
in terms of the non-negative state variables $x$ and $y$ and system parameters $\beta\ge 0$, $\delta,\gamma,\sigma>0$, and with a superscribed dot denoting differentiation with respect to time. These equations model a first-order exothermic chemical reaction in a continuous stirred tank reactor (here, and in the literature, abbreviated as CSTR) with $x$ and $y$ describing the reagent concentration and mixture temperature in the tank, respectively, as described by Bykov and Tsybenova, 2001 \cite{BT01}. 

\begin{figure}[ht]
  \centering
  \includegraphics[width=0.6\textwidth]{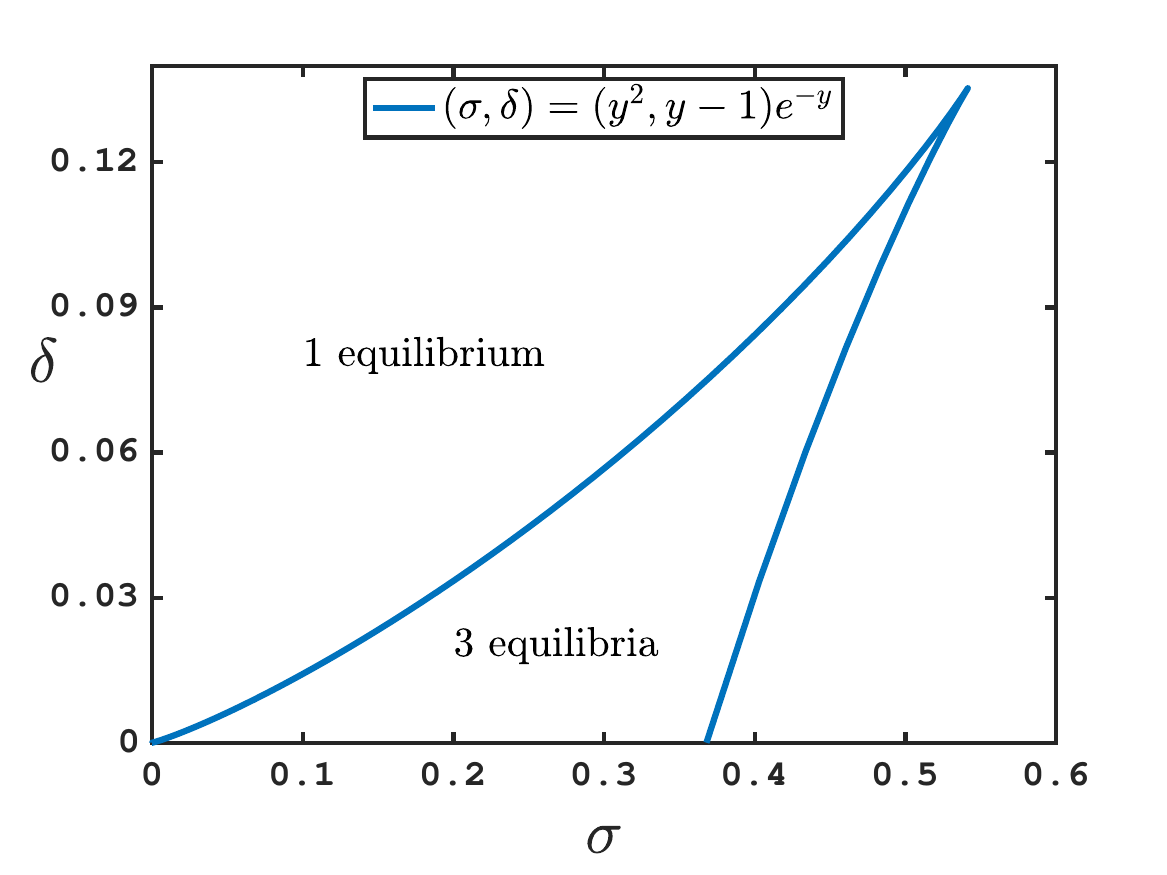}
  \caption[Number of equilibria]{When $\beta=0$, the number of equilibria equals three for $\delta<1/e^2$ and $\sigma\in \left(e^{-y_*}y_*^2, e^{-y^*}y^{*2}\right)$, where $(1-y_*)e^{-y_*}=(1-y^*)e^{-y^*}=-\delta$. Two equilbria are found along the curve $y\mapsto(\sigma,\delta)=\left(y^2,y-1\right)e^{-y}$.}
  \label{fig:cstr:equilibria}
\end{figure}

Equilibria of \eqref{eq:cstr:model} are pairs of values of $x$ and $y$, given allowable values for the parameters $\beta$, $\delta$, $\gamma$, and $\sigma$, for which $\dot{x}$ and $\dot{y}$ both equal $0$. By elimination, a pair $(x,y)$ is an equilibrium provided that $x=\delta y/\sigma$, where
\begin{equation}
\label{eq:cstr:eq1}
    \sigma-\delta y=ye^{-\frac{y}{1+\beta y}},
\end{equation}
i.e., by necessity that $0<y<\sigma/\delta$ and $0<x<1$, since the right-hand side is non-negative. From an analysis of the graph of the right-hand side of \eqref{eq:cstr:eq1} as a function of $y$ (say an exploration of local maxima and minima and limiting behaviors as $y\rightarrow 0$ and $y\rightarrow\infty$) and recognizing the graph of the left-hand side as a straight line with positive ordinate intercept and negative slope, it follows that there exists at least one equilibrium for any allowable choice of $\beta$, $\delta$, $\gamma$, and $\sigma$. The equilibrium is unique provided that either $\beta\ge1/4$, $\delta>(1-4\beta)/e^2$, or
\[
0\le\beta<1/4,\,\delta\le(1-4\beta)/e^2\mbox{, and }\sigma\notin\left[\delta y_*+y_*e^{-\frac{y_*}{1+\beta y_*}},\delta y^*+y^*e^{-\frac{y^*}{1+\beta y^*}}\right]
\]
where $y_*\le y^*$ are the roots of
\[
\frac{1+y^2\beta^2-y(1-2\beta)}{(1+\beta y)^2}e^{-\frac{y}{1+\beta y}}=-\delta.
\]
For all other parameter combinations, there exist at least two and at most three equilibria (see Fig.~\ref{fig:cstr:equilibria} for the case when $\beta=0$). For an analysis of their stability and bifurcations, it is convenient to first restrict attention to a limiting case that permits straightforward closed-form analysis.

\subsection{The case when $\beta=0$}
\label{sec:cstr:beta0}
When $\beta=0$ and for fixed $\gamma$, the collection of equilbria may be described conveniently by the two-dimensional parameterization
\begin{align}
  \label{cstr:eqfamily}
  \left\{\left(x,y,\delta,\sigma\right)\,|\;0<x<1,\delta=\frac{xe^{-y}}{1-x},\,\sigma=\frac{ye^{-y}}{1-x}\right\}.
\end{align}
Their stability is to linear approximation characterized by the eigenvalues of the Jacobian matrix
\begin{equation}
  \label{cstr:Jdef}
  J:=\begin{pmatrix}-e^y/x & e^y(1-x) \\ -e^y/\gamma & -e^y(1-x)(1-y)/\gamma y\end{pmatrix},
\end{equation}
which add to $0$ along a curve with
\begin{equation}
\label{eq:cstr:hb}
y=\frac{x(1-x)}{x(1-x)-\gamma}
\end{equation}
(where $\trace J=0$) and one of which equals $0$ along a curve with
  \[
  y=\frac{1}{1-x}
  \]
(where $\det J=0$).
Provided that $\gamma<4/27$, these curves intersect at distinct points corresponding to two real positive roots of the polynomial $p(x):=x^3-x^2+\gamma$, whereas no such intersections exist for $\gamma>4/27$. We denote these points of intersection by the subscripts $_{\mathrm{BT},1}$ and $_{\mathrm{BT},2}$, respectively. Points along the first curve with $x_{\mathrm{BT},1}<x<x_{\mathrm{BT},2}$ (where $p(x)<0$) correspond to \textbf{Hopf bifurcations} (for which the two eigenvalues of $J$ are both imaginary and each other's negatives) along one-parameter families of equilibria that intersect the curve transversally. Similarly, all points on the second curve are \textbf{saddle-node bifurcations} along one-parameter families of equilibria that intersect the curve transversally. Their intersections at $x_{\mathrm{BT},1}$ and $x_{\mathrm{BT},2}$ are \textbf{Bogdanov-Takens bifurcations} that coincide for $\gamma=4/27$. 

Consistent with the implementation in \textsc{coco}, which also generalizes to problems of higher state-space dimension, we may identify Hopf bifurcations in this model problem as a subset of those equilibria for which
\begin{equation}
\label{eq:cstr:hopf}
J^2\cdot \begin{pmatrix}\cos\theta\\\sin\theta\end{pmatrix}+k \begin{pmatrix}\cos\theta\\\sin\theta\end{pmatrix}=0
\end{equation}
for $k>0$ and independent of $\theta$. Solving \eqref{eq:cstr:hopf} using the expression \eqref{cstr:Jdef} for $J$, we obtain the $\theta$-independent one-dimensional family of solutions given by
\begin{equation}
  \label{cstr:hopfsol}
  k=-\frac{p(x)}{x^2\gamma}e^{2y},\,y=\frac{x(1-x)}{x(1-x)-\gamma}
\end{equation}
with $p(x)$ given above. It follows that $J^2$ is a multiple of the $2\times2$ identity matrix along the curve \eqref{eq:cstr:hb}, specifically that $J^2=-kI_2$ along this curve. 
As expected, $k>0$ for $x\in(x_{\mathrm{BT},1},x_{\mathrm{BT},2})$, such that those solutions correspond to Hopf bifurcations. On the other hand, values of $x$ with $k<0$  correspond to \textbf{neutral saddles} (for which the two eigenvalues of $J$ are both real and sum to $0$). 

System~\eqref{eq:cstr:hopf} also has the $\theta$-dependent two-dimensional family of solutions given by
\begin{align}
  \label{cstr:non-neutral}
  k=-\frac{(1-x(1-x)\tan\theta)^2}{x^2}e^{2y},\,y=\left(1+\frac{\gamma}{x(1-x)}-\frac{\cot\theta}{1-x}-\gamma\tan\theta\right)^{-1},
\end{align}
along which $k$ is always non-positive. Along this family, $-\sqrt{-k}$ is an eigenvalue of $J$ with eigenvector $(\cos\theta,\sin\theta)^\mathsf{T}$ while the other eigenvalue equals $e^y(x\cot\theta-\gamma)/x\gamma$. For a given $\theta$, the two families in \eqref{cstr:hopfsol} and \eqref{cstr:non-neutral} intersect on a neutral saddle with
\[
x=x_{\mathrm{BP},\pm}(\theta):=\frac{1}{2\gamma}\left(\gamma+\cot^2\theta\pm\sqrt{(\gamma+\cot^2\theta))^2-8\gamma^2\cot\theta}\right).
\]
For later reference, we anticipate the possible detection of such \textbf{branch points} during analysis with \textsc{coco}.

\begin{figure}[ht]
  \centering
  \includegraphics[width=0.6\textwidth]{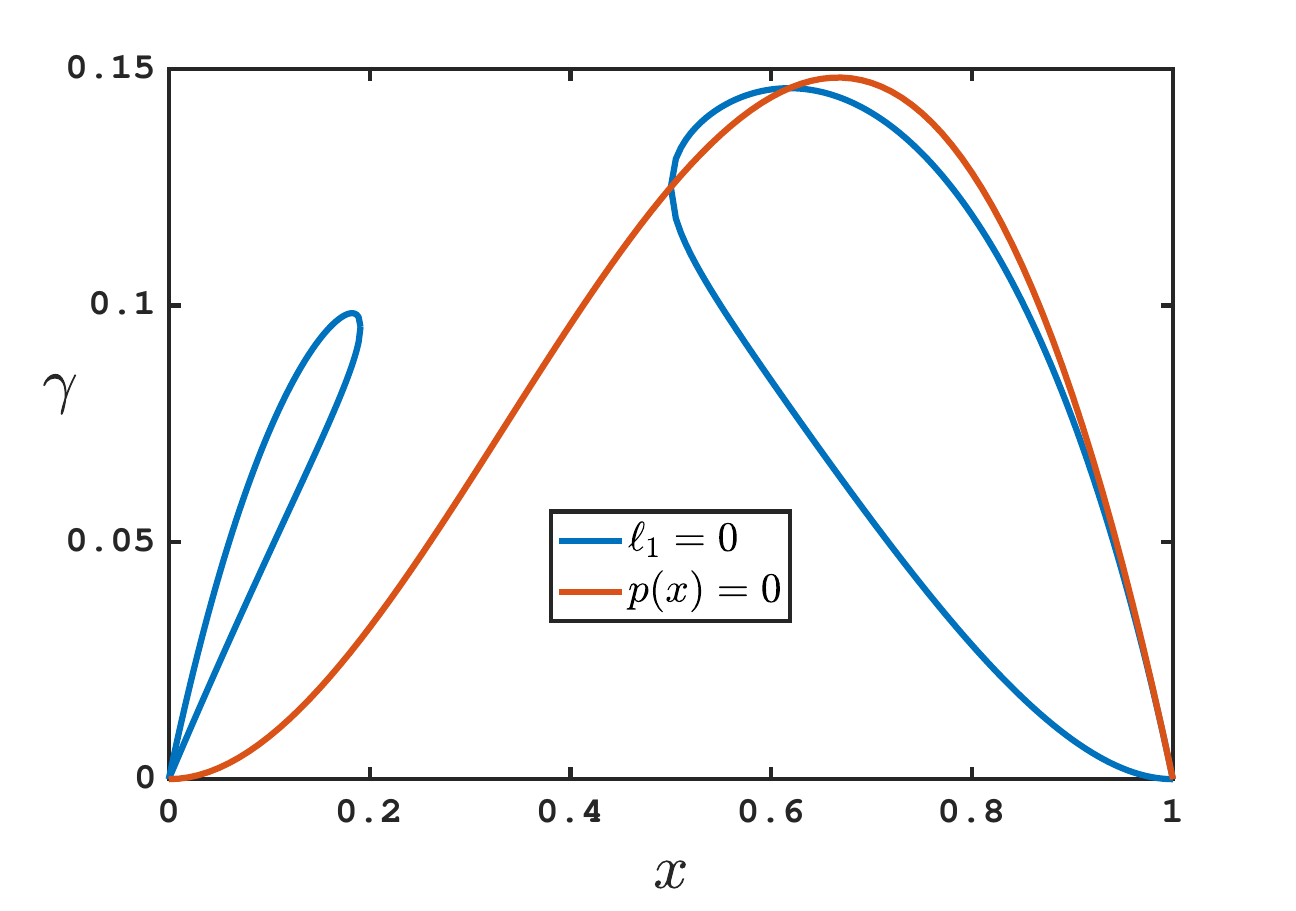}
  \caption[Loci of degenerate Hopf bifurcations]{For $x_{\mathrm{BT},1}<x<x_{\mathrm{BT},2}\Leftrightarrow p(x)<0$, the $\ell_1=0$ level set consists of degenerate Hopf bifurcations. In the figure, the $p(x)=0$ level set consists of families of Bogdanov-Takens points that coincide for $\gamma=4/27$. Supercritical Hopf bifurcations occur for $(x,\gamma)$ inside the right-most region enclosed by the $\ell_1=0$ level set and below $p(x)=0$.}
  \label{fig:cstr:lyapunov}
\end{figure}

Along the subset of Hopf bifurcations, following lengthy derivations, the \textbf{first Lyapunov coefficient} evaluates to
\[
\ell_1:=\frac{2 \gamma^2-x(3-\gamma)\gamma+x^2(1+5\gamma-2\gamma^2)-x^3(3+2\gamma)+3x^4-x^5}{4\sqrt{\gamma}(x^2(1-x)-\gamma)^{3/2}(1+(1-x)\gamma)}.
\]
With reference to Fig.~\ref{fig:cstr:lyapunov}, by solving the equation $\ell_1=0$ for $\gamma$ while recalling that $\gamma<x^2(1-x)$ (since $p(x)<0$), one concludes that the first Lyapunov coefficient vanishes for $x_{\mathrm{BT},1}<x<x_{\mathrm{BT},2}$ at a unique value of $x$ when $1/8\le\gamma<(7-3\sqrt{5})/2$, at two distinct values of $x$ when $0<\gamma<1/8$, and for no value of $x$ when $\gamma\ge(7-3\sqrt{5})/2$. Such loci of vanishing first Lyapunov coefficient, denoted by the subscript $_\mathrm{DH}$, correspond to \textbf{degenerate Hopf bifurcations}, such that the Hopf bifurcations on either side are \textbf{supercritical} where $\ell_1<0$ and \textbf{subcritical} where $\ell_1>0$. As $\gamma$ approaches $(7-3\sqrt{5})/2$ from below, the single $x_\mathrm{DH}\rightarrow x_{\mathrm{BT},1}$. Similarly, the left-most $x_{\mathrm{DH}}\rightarrow x_{\mathrm{BT},1}$ as $\gamma$ approaches $1/8$ from below. Away from these degenerate values of $\gamma$, it is evident from the closed-form expression that $\left|\ell_1\right|$ grows without bound as $x$ approaches $x_{\mathrm{BT},1}$ and $x_{\mathrm{BT},2}$.

\subsection{Numerical bifurcation analysis}
\label{sec:cstr:num}
The analytical results of the previous section may be explored numerically and then extended to $\beta\ne0$ using a problem encoding in \textsc{coco}. This may be found in the previously referenced code repository. The text below includes extracts of the associated screen output with accompanying commentary.

\begin{figure}[ht]
  \centering
  \includegraphics[width=1\textwidth]{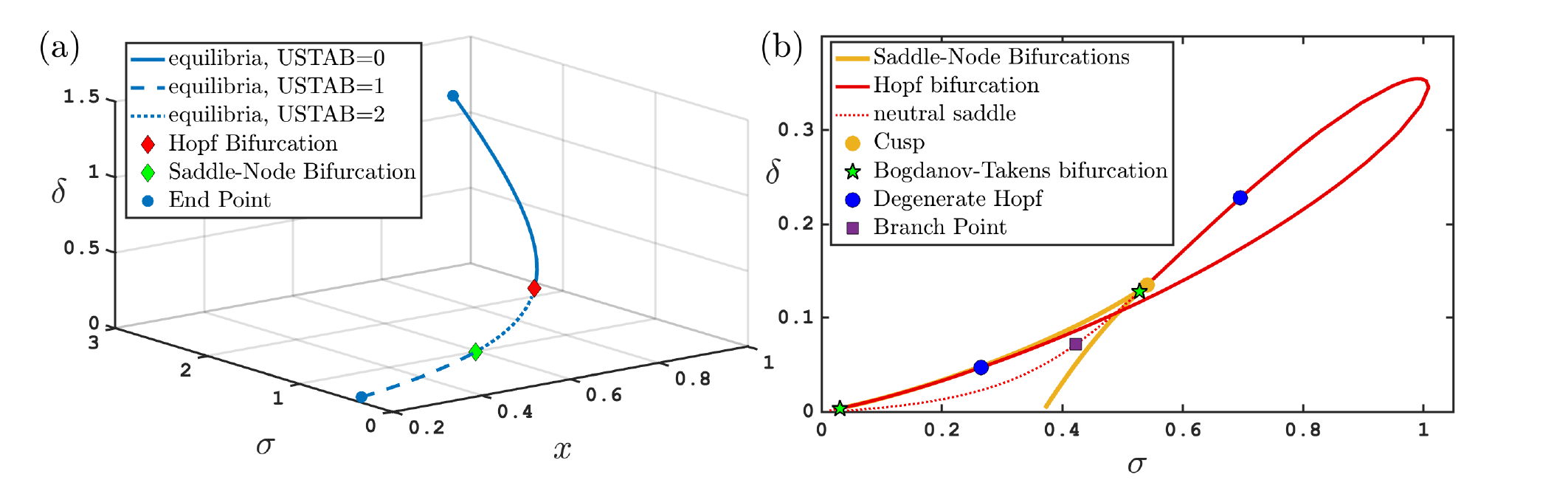}
  \caption[Branches of equilibria in CSTR]{(a) Branch of equilibria in
    the CSTR problem \eqref{eq:cstr:model}, obtained by fixing $y$ while
    varying $\delta$ and $\sigma$. Line style indicates
    number of unstable eigenvalues (\texttt{USTAB}). Hopf (red diamond)
    and saddle-node (green diamond) bifurcations are indicated along the
    branch. (b) Curves of saddle-node and Hopf bifurcations of equilibria in the
    $(\sigma,\delta)$-plane. Along the curves cusp, Bogdanov-Takens,
    and degenerate Hopf bifurcation points are indicated. A curve where the
    saddle is neutral starts at the Bogdanov-Takens points and
    contains a branch point of the defining problem. Other parameters: $\beta=0$,
    $\gamma=0.1$.}
  \label{fig:equilibria}
\end{figure}

To begin with, the screen output below shows the result of numerical continuation along a branch of equilibria for constant values of $y$, $\beta$, and $\gamma$, and with initial solution guess given by the equilibrium at $(x,y)=(9/10,2)$ obtained for $(\beta,\gamma,\delta,\sigma)=(0,1/10,9/e^2,20/e^2)$.
\begin{lstlisting}[language=coco-highlight,frame=lines]
...   LABEL  TYPE             x        sigma        delta
...       1  EP      9.0000e-01   2.7067e+00   1.2180e+00
...       2  HB      7.2361e-01   9.7930e-01   3.5431e-01
...       3  SN      5.0000e-01   5.4134e-01   1.3534e-01
...       4  EP      2.0000e-01   3.3834e-01   3.3834e-02
\end{lstlisting}
The branch is visualized in Fig.~\ref{fig:equilibria}(a) with line style indicating the number of eigenvalues with positive real part. The saddle-node and Hopf bifurcations denoted by \mcode{SN} and \mcode{HB} in the screen output agree to the number of digits shown with the theoretical predictions
\[
(x_\mathrm{SN},\delta_\mathrm{SN},\sigma_\mathrm{SN})=\left(\frac{1}{2},\frac{1}{e^2},\frac{4}{e^2}\right)
\]
and
\[
(x_\mathrm{HB},\delta_\mathrm{HB},\sigma_\mathrm{HB})=\left(\frac{5+\sqrt{5}}{10},\frac{3+\sqrt{5}}{2e^2},\frac{5+\sqrt{5}}{e^2}\right).
\]
Notably, since $y$ is held fixed while $\sigma$ and $\delta$ both vary, the saddle-node bifurcation does not coincide with a local maximum in either of the two problem parameters, as would have been expected if $y$ were free to vary and either $\sigma$ or $\delta$ were held fixed\footnote{We take this opportunity to point out a ``bug'' in the encoding of earlier releases of \textsc{coco}. This assumed a relationship between the tangent vector to the solution manifold at a saddle-node bifurcation point and the eigenvector of the Jacobian of the vector field corresponding to the zero eigenvalue, as would be the case if the state variables and only one of the system parameters were free to vary. This bug affected none of the demos included in the earlier release, but is triggered precisely by the first continuation run in this section which holds $y$ fixed. The bug has been corrected in the release of \textsc{coco} accompanying this chapter.}.

Continuation along branches of saddle-node and Hopf bifurcations, respectively, intersecting the points found above and allowing variations also in $y$, yields the curves shown in Fig.~\ref{fig:equilibria}(b) in a projection onto the $(\sigma,\delta)$ parameter plane. In particular, the saddle-node bifurcation curve agrees with the curve shown in Fig.~\ref{fig:cstr:equilibria}, while the Hopf bifurcation curve agrees with the theoretical prediction
\[
[0,1)\ni s\mapsto\left(\frac{se^{-\frac{s(1-s)}{s(1-s)-\gamma}}}{s(1-s)-\gamma},\frac{se^{-\frac{s(1-s)}{s(1-s)-\gamma}}}{1-s}\right)
\]
from Section~\ref{sec:cstr:beta0}. The screen output during continuation along the branch of Hopf bifurcations in Fig.~\ref{fig:equilibria}(b) is shown below.
\begin{lstlisting}[language=coco-highlight,frame=lines]
...   LABEL  TYPE         sigma        delta            x            y
...       1  EP      9.7930e-01   3.5431e-01   7.2361e-01   2.0000e+00
...       2  DH      6.9586e-01   2.2796e-01   5.4964e-01   1.6778e+00
...       3  BTP     5.2818e-01   1.2801e-01   4.1261e-01   1.7024e+00
...       4  BP      4.2292e-01   7.1751e-02   3.1582e-01   1.8615e+00
...       5  EP      2.3161e-01   1.7371e-02   2.0000e-01   2.6667e+00

...   LABEL  TYPE         sigma        delta            x            y
...       6  EP      9.7930e-01   3.5431e-01   7.2361e-01   2.0000e+00
...       7  FP      1.0082e+00   3.4657e-01   7.5943e-01   2.2092e+00
...       8  DH      2.6535e-01   4.7034e-02   8.5174e-01   4.8052e+00
...       9  BTP     3.0747e-02   3.5466e-03   8.6695e-01   7.5160e+00
...      10  EP      1.6578e-02   1.7480e-03   8.6904e-01   8.2418e+00
\end{lstlisting}
Here, the two Bogdanov-Takens points denoted by \mcode{BTP} agree to the number of digits shown with the theoretical prediction obtained from the two positive roots of the polynomial $p(x)$ defined in the previous section. Similarly, the two degenerate Hopf bifurcations denoted by \mcode{DH} agree to the number of digits shown with the theoretical prediction obtained from the two roots of $\ell_1=0$ in $(0,1)$.

The point denoted by \mcode{FP} is a local maximum in the value of $\sigma$, while the branch point denoted by \mcode{BP} represents a point of intersection along the family of neutral saddles with a family of solutions to \eqref{eq:cstr:hopf} for some value of $\theta$. Indeed, to the number of digits shown, this agrees with the predicted value $x_\mathrm{BP,-}(\theta)$ for $\theta\approx 1.43$ obtained from the vector $(\cos\theta,\sin\theta)^\mathsf{T}$ stored by \textsc{coco} at this point. 

Allowing variations also in $\gamma$, a final stage of continuation along a branch of degenerate Hopf bifurcations starting from from $(\sigma,\delta)\approx (0.70,0.23)$ shows such bifurcations persisting until $\gamma$ gets close to the predicted value of $1/8$ to the number of digits shown in the screen output below.
\begin{lstlisting}[language=coco-highlight,frame=lines]
...   LABEL  TYPE         sigma        delta        gamma
...       1  EP      6.9586e-01   2.2796e-01   1.0000e-01
...       2  MX      5.4135e-01   1.3534e-01   1.2500e-01

...   LABEL  TYPE         sigma        delta        gamma
...       3  EP      6.9586e-01   2.2796e-01   1.0000e-01
...       4  EP      7.9678e-01   3.0000e-01   8.7404e-02
\end{lstlisting}
The continuation algorithm terminates with a failure to converge (denoted by \mcode{MX}) as, near $\gamma=1/8$, $\ell_1$ varies rapidly between $0$ at the \mcode{DH} point and $\infty$ at the \mcode{BTP} point (which are coincident in the limit of $\gamma=1/8$).

Similar observations apply to the branch of degenerate Hopf bifurcations starting from $(\sigma,\delta)\approx(0.26,0.047)$, which persists until $\gamma$ gets close to $(7-3\sqrt{5})/2$ to the number of digits shown in the screen outputs below.
\begin{lstlisting}[language=coco-highlight,frame=lines]
...   LABEL  TYPE         sigma        delta        gamma
...       1  EP      2.6535e-01   4.7034e-02   1.0000e-01
...       2  EP      2.1043e-01   3.4844e-02   9.2368e-02

...   LABEL  TYPE         sigma        delta        gamma
...       3  EP      2.6535e-01   4.7034e-02   1.0000e-01
...       4  FP      4.9998e-01   1.1803e-01   1.4590e-01
...       5  MX      4.9998e-01   1.1803e-01   1.4590e-01
\end{lstlisting}
The spurious fold point observed here is also a consequence of the nondifferentiability of $\ell_1$ in the limit of $\gamma=(7-3\sqrt{5})/2$.

\subsection{Periodic orbits}
\label{sec:cstr:per}
Branches of periodic orbits emanate from points on the Hopf bifurcation curve found above under variations in $\sigma$ and fixed $\beta,\gamma,\delta$. These curves either form \textbf{Hopf bubbles} that connect two distinct Hopf bifurcation points or terminate on a \textbf{homoclinic orbit} that limits on the same equilibrium in both forward and backward time. In each case, their numerical study is initiated from an initial solution guess that deviates from the Hopf bifurcation equilibrium by a small-radius circular orbit of the form $v\cos\omega t-w\sin\omega t$ for a unique $\omega$ and any vectors $v$ and $w$ satisfying $Jv=-\omega w$ and $Jw=\omega v$ for $J$ evaluated at the Hopf bifurcation.

\begin{figure}[ht]
  \centering
  \includegraphics[width=0.85\textwidth]{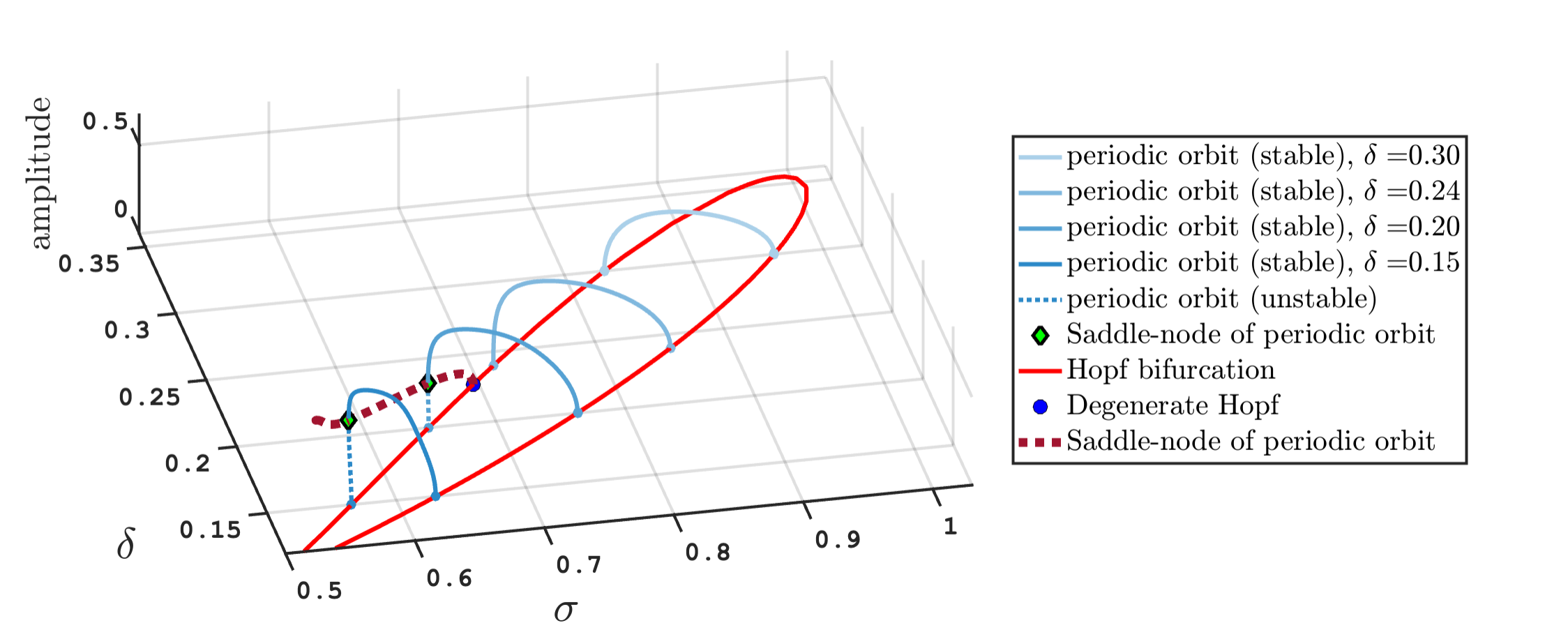}
  \caption{Branches of periodic orbits connecting Hopf bifurcations (Hopf bubbles) under variations in $\sigma$ and for different values of $\delta\geq0.15$. For $\delta$ below the value associated with the degenerate Hopf bifurcation, such bubbles include a saddle-node bifurcation of periodic orbits. Other parameters: $\beta=0$, $\gamma=0.1$.}
  \label{fig:hopfbubbles}
\end{figure}

The screen output below shows the result of numerical continuation along a branch of periodic orbits emanating from the supercritical Hopf bifurcation found at $(\sigma,\delta,x,y)\approx(0.85, 0.24, 0.81, 2.90)$ where $2\pi/\omega\approx0.57$.
\begin{lstlisting}[language=coco-highlight,frame=lines]
...   LABEL  TYPE         sigma    po.period    amplitude
...       1  EP      8.5000e-01   5.7480e-01   2.7993e-04
...       2          8.3465e-01   6.1989e-01   2.0863e-01
...       3          7.9414e-01   7.3451e-01   3.8426e-01
...       4          7.6497e-01   8.0974e-01   4.4289e-01
...       5          7.3772e-01   8.9386e-01   4.5054e-01
...       6          7.1521e-01   1.0332e+00   2.4601e-01
...       7  EP      7.1379e-01   1.0583e+00   6.6074e-03
\end{lstlisting}
Here, the column denoted by \mcode{po.period} contains the period of each orbit, while the column denoted by \mcode{amplitude} contains the difference between the maximum and minimum values of $x$ for each orbit. As shown in Fig.~\ref{fig:hopfbubbles}, this branch of periodic orbits terminates on a supercritical Hopf bifurcation found at $(\sigma,\delta,x,y)\approx(0.71, 0.24, 0.56, 1.68)$ where $2\pi/\omega\approx1.06$.

The screen output below shows the result of numerical continuation along a branch of periodic orbits emanating from the supercritical Hopf bifurcation found at $(\sigma,\delta,x,y)\approx(0.63, 0.15, 0.83, 3.48)$ where $2\pi/\omega\approx0.39$.
\begin{lstlisting}[language=coco-highlight,frame=lines]
...   LABEL  TYPE         sigma    po.period    amplitude
...       1  EP      6.3433e-01   3.9398e-01   2.7389e-04
...       2          6.0924e-01   7.0868e-01   5.1197e-01
...       3          5.7782e-01   1.0020e+00   6.3901e-01
...       4  SN      5.6537e-01   1.4793e+00   4.7242e-01
...       5          5.6588e-01   1.5938e+00   3.1561e-01
...       6  EP      5.6748e-01   1.5989e+00   9.3532e-03
\end{lstlisting}
As shown in Fig.~\ref{fig:hopfbubbles}, after undergoing a saddle-node bifurcation denoted by \mcode{SN}, this branch of periodic orbits terminates on a subcritical Hopf bifurcation found at $(\sigma,\delta,x,y)\approx(0.57,0.15, 0.45, 1.68)$ where $2\pi/\omega\approx1.60$.

As illustrated in Fig.~\ref{fig:hopfbubbles} below, the transition between Hopf bubbles connecting supercritical Hopf bifurcations and those connecting super- and subcritical Hopf bifurcations occurs at the degenerate Hopf bifurcation found at $\delta\approx0.23$. This bifurcation is also a terminal point for the branch of saddle-node bifurcations of periodic orbits shown in the figure.

\begin{figure}[ht]
  \centering
  \includegraphics[width=\textwidth]{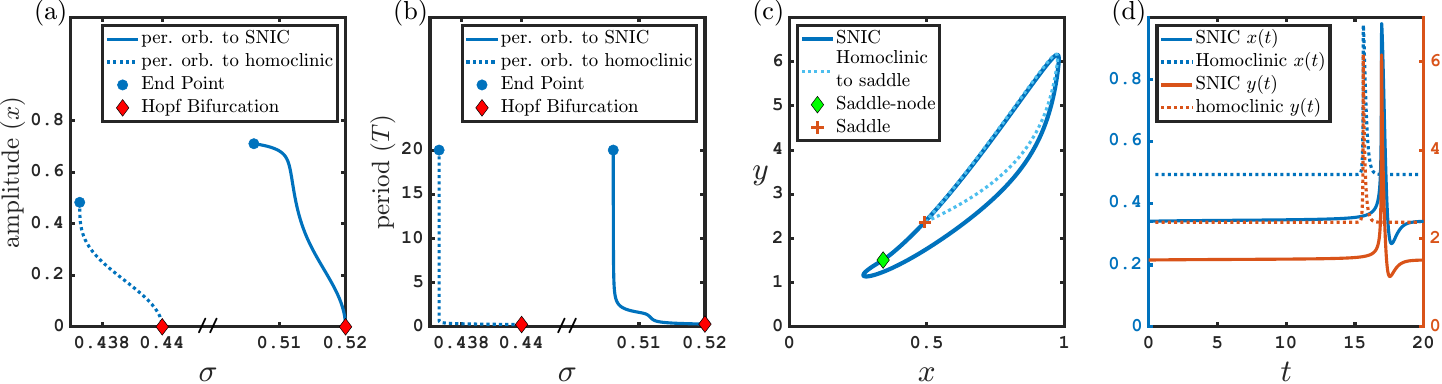}
  \caption{Branches of periodic orbits emanating from Hopf bifurcations at
    $(\sigma,\delta)=(0.52,0.11)$ (shown in Fig.~\ref{fig:equilibria}(a)) and $(\sigma,\delta)=(0.44,0.091)$,
    for fixed $\delta$ and varying $\sigma$. (a,b) Bifurcation diagrams
    in $(\sigma,\max(x)-\min(x))$-plane and $(\sigma,T)$-plane, where $T$ denotes the period. (c,d)
    Phase portrait and time profiles of periodic orbits with period $20$ close to the homoclinic and SNIC bifurcation, respectively. Other parameters:
    $\beta=0$, $\gamma=0.1$.}
  \label{fig:homsnic}
\end{figure}

For sufficiently small values of $\delta$, the Hopf bubbles rupture to form two independent branches of periodic orbits. The screen output below shows the result of numerical continuation along a branch of periodic orbits emanating from the supercritical Hopf bifurcation found at $(\sigma,\delta,x,y)\approx(0.52,0.11,0.84,3.81)$ where $2\pi/\omega\approx3.15$.
\begin{lstlisting}[language=coco-highlight,frame=lines]
...   LABEL  TYPE     po.period        sigma    amplitude
...       1  EP      3.1454e-01   5.1999e-01   2.7407e-04
...       2  UZ      5.0000e+00   5.0623e-01   7.1006e-01
...       3  UZ      1.0000e+01   5.0614e-01   7.1039e-01
...       4  UZ      1.5000e+01   5.0613e-01   7.1044e-01
...       5  EP      2.0000e+01   5.0612e-01   7.1040e-01
\end{lstlisting}
As shown in Fig.~\ref{fig:homsnic}, this branch of periodic orbits terminates on a homoclinic orbit to a saddle-node bifurcation equilibrium at $(\sigma,\delta,x,y)\approx(0.51,0.11, 0.34, 1.53)$, corresponding to a so-called \textbf{saddle-node on invariant circle (SNIC) bifurcation}.  To make visible the rapid increase of the period over a small parameter range near the SNIC bifurcation, the screen output reports solutions where the period is a multiple of $5$ (denoted by \mcode{UZ}). 

Finally, the screen output below shows the result of numerical continuation along a branch of periodic orbits emanating from the supercritical Hopf bifurcation found at $(\sigma,\delta,x,y)\approx(0.44, 0.091,0.84, 4.07)$ where $2\pi/\omega\approx0.26$.
\begin{lstlisting}[language=coco-highlight,frame=lines]
...   LABEL  TYPE     po.period        sigma    amplitude
...       1  EP      2.6304e-01   4.4000e-01   2.7567e-04
...       2  UZ      5.0000e+00   4.3724e-01   4.8302e-01
...       3  UZ      1.0000e+01   4.3724e-01   4.8302e-01
...       4  UZ      1.5000e+01   4.3724e-01   4.8302e-01
...       5  EP      2.0000e+01   4.3724e-01   4.8302e-01
\end{lstlisting}
As shown in Fig.~\ref{fig:homsnic}, this branch of periodic orbits terminates on a homoclinic orbit to a saddle equilibrium at a \textbf{homoclinic bifurcation} for $\sigma\approx 0.437$. The period appears to increase more rapidly than when approaching the SNIC as the parameter range in $\sigma$ is noticeably smaller.

Continuation along a branch of homoclinic orbits under simultaneous variations in $\sigma$ and $\delta$ may be accomplished using a periodic orbit approximant of very high period. Such a periodic orbit will pass near either a saddle or a saddle-node equilibrium, such that it is instructive to monitor the determinant and trace of the Jacobian of the vector field $(f_x,f_y)$ in \eqref{eq:cstr:model} at the point of the orbit with smallest value of $f_x^2+f_y^2$.

The screen output below shows the results of such continuation of an orbit with period $500$ with initial solution guess constructed from the terminal orbit in the third run in this section.
\begin{lstlisting}[language=coco-highlight,frame=lines]
... LABEL  TYPE       sigma       delta   amplitude          det           tr
...     1  EP    5.0612e-01  1.1426e-01  7.1045e-01   6.7738e-04  -2.9728e+00
...     2  NCS   5.0027e-01  1.1038e-01  6.5464e-01  -5.0000e-02  -3.7011e+00
...     3  NSA   4.9376e-01  1.0821e-01  6.0200e-01  -2.2875e+01   2.9426e-07
...     4  MX    3.0747e-02  3.5466e-03  1.3918e-07  -2.0245e-01   7.4677e-04

... LABEL  TYPE       sigma       delta   amplitude          det           tr
...     5  EP    5.0612e-01  1.1426e-01  7.1045e-01   6.7738e-04  -2.9728e+00
...     6  NCS   5.0918e-01  1.1625e-01  7.0956e-01  -5.0000e-02  -2.5766e+00
...     7  NSA   5.1949e-01  1.2293e-01  6.9552e-01  -4.7361e+00   1.1280e-07
...     8  NCS   5.3240e-01  1.3122e-01  5.9868e-01  -5.0000e-02   6.5372e+00
...     9  NCS   5.3550e-01  1.3261e-01  1.9258e-01  -5.4838e-02   5.9755e+00
...    10  FP    5.3561e-01  1.3264e-01  1.7202e-01  -1.5543e+00   5.3676e+00
...    11  NCS   5.2825e-01  1.2805e-01  8.2863e-04  -5.0000e-02   2.6269e-02
...    12  NSA   5.2818e-01  1.2801e-01  3.3022e-06   4.8011e-04  -8.6704e-05
...    13  BP    5.2818e-01  1.2801e-01  5.5482e-07   2.1720e-04  -2.8064e-04
...    14  MX    5.2816e-01  1.2800e-01  2.8784e-10   1.5946e-04  -3.3576e-03
\end{lstlisting}
The computation terminates with failure of convergence (\mcode{MX}) at the two Bogdanov-Takens points, identifiable by the small amplitude of the orbit.

\begin{figure}[ht]
  \centering
  \includegraphics[width=1\textwidth]{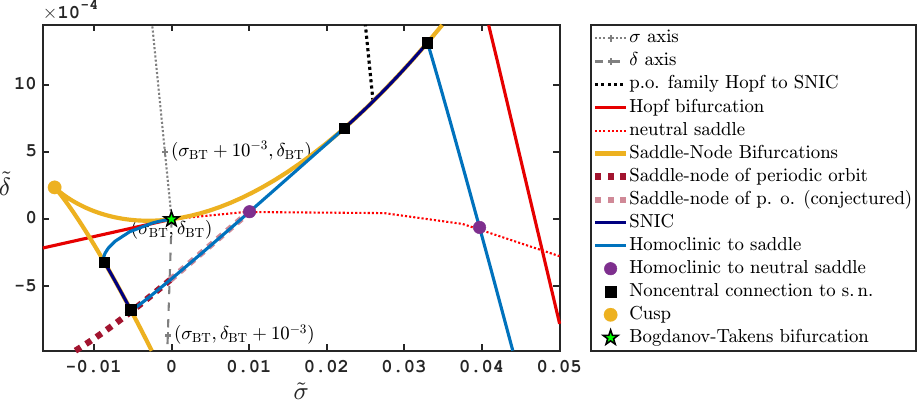}
  \caption[Bifurcation diagram for equilibria and periodic orbits in
  two parameters]{Bifurcation diagram for equilibria and periodic orbits with free parameters $\sigma$ and $\delta$ in a small parameter region containing the Bogdanov-Takens point $(\sigma_\mathrm{BT,1},\delta_\mathrm{BT,1})$ (subscript $1$ not shown in the figure) in  Fig.~\ref{fig:equilibria}(b). For better visibility of the bifurcation curves, the coordinates of the parameter plane have been rotated counterclockwise by $\theta=150^\circ$ around the Bogdanov-Takens point, such that $\tilde{\sigma}=\cos\theta(\sigma-\sigma_\mathrm{BT})-\sin\theta(\delta-\delta_\mathrm{BT})$,    $\tilde{\delta}=\sin\theta(\sigma-\sigma_\mathrm{BT})+\cos\theta(\delta-\delta_\mathrm{BT})$. Other parameters: $\beta=0$, $\gamma=0.1$.}
  \label{fig:cstr:bif2d}
\end{figure}

Changes to the sign of the trace (\mcode{tr}), here identified by the label \mcode{NSA}, correspond approximately to homoclinic orbits connecting to neutral saddle equilibria. Fig.~\ref{fig:cstr:bif2d} shows two such points in the vicinity of the Bogdanov-Takens point $(\sigma_\mathrm{BT,1},\delta_\mathrm{BT,1})$ at points of intersection of the branch of neutral saddles and the branch of homoclinic orbits.

For homoclinic orbits connecting to saddle-node bifurcation equilibria, the determinant (\mcode{det}) should equal $0$, whereas the determinant is negative in the saddle case. Across the transition between such orbits, the shape of the homoclinic turns from a smooth curve to a non-smooth curve in phase space (the ``corner'' is the location of the saddle, see phase portraits in Fig.~\ref{fig:homsnic}). At the transition point, the homoclinic orbit is a \textbf{non-central connection to a saddle-node} \cite{guckenheimer1986multiple}. In the screen output, approximate detection of such a transition is accomplished by monitoring for a crossing by \mcode{det} of $-0.05$ (a suitably chosen small, negative number) and the corresponding points are denoted by \mcode{NCS}. Four such non-central connections to saddle-nodes are found near the Bogdanov-Takens point $(\sigma_\mathrm{BT,1},\delta_\mathrm{BT,1})$ at the end points of the segments labeled SNIC in Fig.~\ref{fig:cstr:bif2d}. Magnification (not included here but available using the provided code) shows that the homoclinic curve is smooth across these points. 

The branch of saddle-node bifurcations of periodic orbits, shown in Fig.~\ref{fig:hopfbubbles} emanating from the degenerate Hopf bifurcation, may be continued to the vicinity of the Bogdanov-Takens point $(\sigma_\mathrm{BT,1},\delta_\mathrm{BT,1})$. Indeed, as shown in Fig.~\ref{fig:cstr:bif2d}, this branch becomes tangential to the curve of homoclinic bifurcations to a saddle equilibrium and is expected (per generic unfolding theory) to terminate at the homoclinic orbit to a neutral saddle equilibrium, even though the numerical algorithm finds it increasingly difficult to converge near this point and terminates prematurely.

\subsection{Reflections and Outlook}
\label{sec:cstr:outlook}
A diagram such as Fig.~\ref{fig:cstr:bif2d} maps out a skeleton of boundaries between different qualitative behaviors of a dynamical system in a two-dimensional parameter plane. Here, variations of $\sigma$ and $\delta$ (or the rotated versions $\tilde{\sigma}$ and $\tilde{\delta}$) along arbitrary curves through this diagram that cross such boundaries result in the appearance, disappearance, or change in stability of equilibria and/or periodic orbits. Further numerical exploration may enable an analyst to build a comprehensive understanding of the possible dynamics of the CSTR model for wide ranges of problem parameters.

\begin{figure}[ht]
  \centering
  \includegraphics[width=1\textwidth]{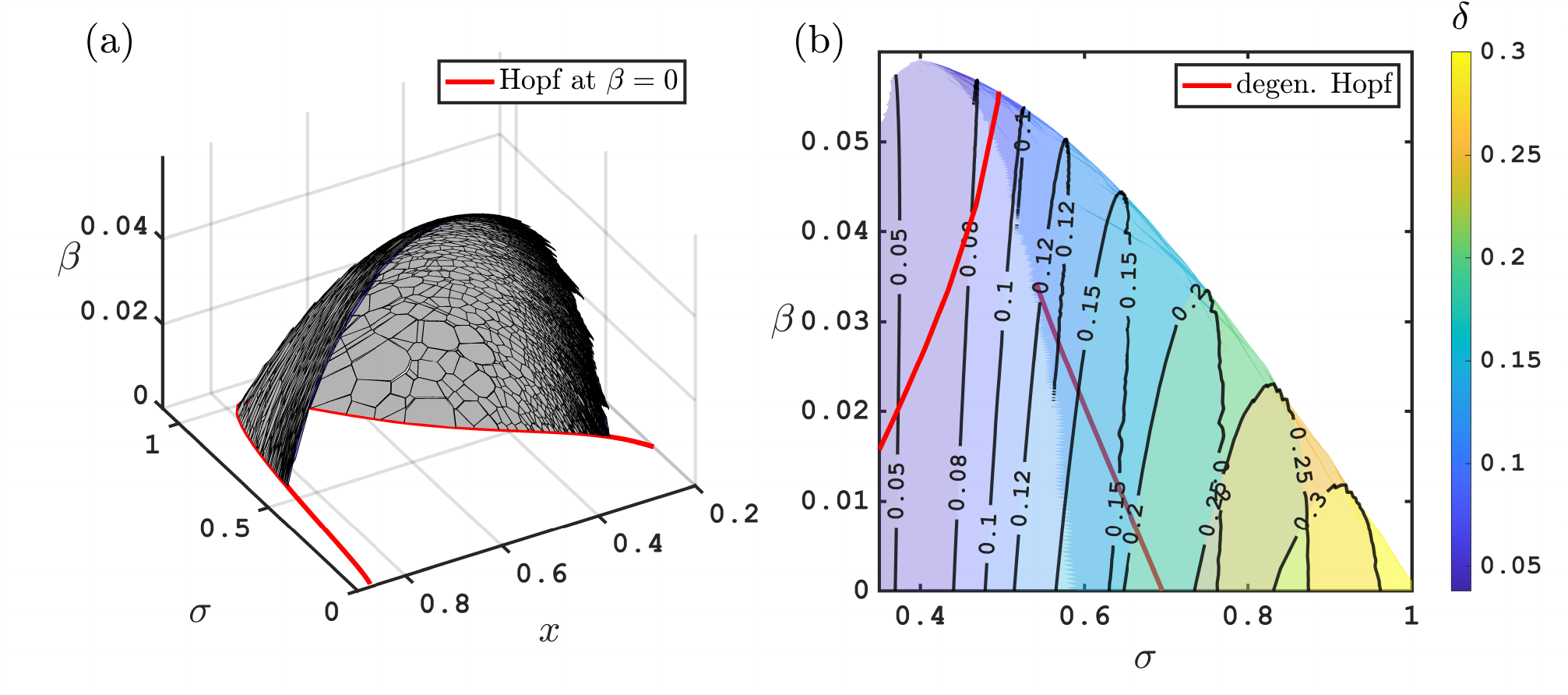}
  \caption{Surface of Hopf bifurcations. (a) Adaptive polygonal tiling
    produced by the manifold growing algorithm described in \cite{dankowicz2020multidimensional,henderson2002multiple}, projected onto $(x,\sigma,\beta)$-space. The red curve in the $\beta=0$ plane corresponds to the curve of Hopf bifurcations in Fig.~\ref{fig:equilibria}(b). (b) Projection onto $(\sigma,\beta,\delta)$-space using contours and including curves of degenerate Hopf bifurcations (starting from points detected in Fig.~\ref{fig:equilibria}(b)). Other parameters: $\gamma=0.1$.}
  \label{fig:cstr:hopf2}
\end{figure}

The skeleton of boundaries obtained for $\beta=0$ and $\gamma=0.1$ in Fig.~\ref{fig:cstr:bif2d} will naturally change under variations in either $\beta$ or $\gamma$. For example, it is anticipated that there exists a combination of parameter values, for which the cusp and Bogdanov-Takens bifurcation point coalesce into a single point, a \textbf{Bogdanov-Takens cusp interaction} \cite{Dumortier_Roussarie_Sotomayor_1987}. The bifurcation diagram in Fig.~\ref{fig:cstr:bif2d} is a two-parameter cross section near this more complex codimension-$3$ bifurcation. Continuation may be used to trace out branches of any of the singular points (the cusp, the homoclinic to a neutral saddle, the non-central connection to a saddle node, and so on) shown in Fig.~\ref{fig:cstr:bif2d} under simultaneous variations in three parameters, as was done for the branch of degenerate Hopf bifurcations in the last runs in Section~\ref{sec:cstr:num}. Similarly, continuation may be used to trace out two-dimensional families of the solution types along any of the curves in Fig.~\ref{fig:cstr:bif2d} under simultaneous variations in three parameters. As an example, panel (a) of Fig.~\ref{fig:cstr:hopf2} shows a two-dimensional surface of Hopf bifurcations under simultaneous variations of $\sigma$, $\delta$, and $\beta$ and projected onto to the $(\sigma,x,\beta)$ space. In panel (b), the same surface is shown in a projection onto the $(\sigma,\beta)$ plane with color and level curves for $\delta$ and overlaid curves of degenerate Hopf bifurcations. 

By now, the reader's impression of computational bifurcation analysis is hopefully of an algorithm that proceeds in stages, starting from facts that may be available through closed-form analysis or proven under suitable assumptions of genericity, and expanding into various directions of inquiry as prompted by observations made in each stage. In this respect, the stages of analysis employed in this section for the CSTR model and, to a degree, even the qualitative nature of the results, are emblematic of what might be expected across a large class of problems.

\section{Defining systems, initialization, and adaptation}
\label{sec:defining_systems}
Underlying the approach to bifurcation analysis illustrated in the preceding section is a general methodology that integrates steps of construction, initialization, adaptation, exploration, and interpretation. This methodology is informed by a relevant mathematical theory that establishes conditions under which various objects may be identified and tracked and allows the analyst to draw appropriate conclusions from the results of a computation. Such theory also alerts the analyst to non-generic instances where one's honed intuition may fail.

At a high level, the general methodology may be understood in terms of three distinct tasks. First, a defining system of algebraic equations or boundary-value problems must be constructed, initialized, and (as necessary) adaptively updated during the analysis. Second, in the case of computational analysis, an algorithm must be implemented that expands a successively growing family of discrete solutions in directions of interest to the analyst and with appropriate resolution along the family. Third, an algorithm must be implemented that monitors for the occurrence of special solution points and that locates these within some desired accuracy. Computational tools such as \textsc{coco}, \textsc{auto}, and \textsc{matcont} automate all three tasks for commonly encountered scenarios and may also allow further user-defined development to support analysis of new classes of problems. A comprehensive overview of defining systems used in \textsc{matcont} \cite{DGK03} and \textsc{auto} \cite{DCFKSW98} is given in the review by Meijer \emph{et al}.~\cite{MDO09}. The defining systems used by \textsc{coco} are described in the tutorial documentation included with this package. The book \textit{Recipes for Continuation} \cite{DS13} describes a general paradigm for how to implement new classes of defining systems in \textsc{coco}.

\subsection{Algebraic problems}
As an example, for continuation along families of equilibria, solutions are obtained from the defining algebraic problem 
\[
\mathbb{R}^{n_x+n_p}\ni\begin{pmatrix}
    x\\p
\end{pmatrix}\mapsto 0=f(x,p)\in\mathbb{R}^{n_x}
\]
in terms of the \textbf{vector field} $f$. This problem requires no adaptive updates, since neither the form of the equations, the number of unknowns, nor their meaning change during the analysis. Initialization of this problem requires, at least, an initial solution guess $(x_0,p_0)$ and may also include information about a preferred initial tangent vector $(x'(s), p'(s))$ along a branch $(x(s),p(s))$ of such solutions. For the CSTR model in Section~\ref{sec:cstr}, $n_x=2$ and $n_p=4$. \textit{A priori}, one thus expects that every equilibrium point lies on a locally four-dimensional manifold of equilibria. Attention may be restricted to an embedded submanifold of lower dimension, for example, by fixing components of $x$ and/or $p$. As a case in point, the two-dimensional family in \eqref{cstr:eqfamily} was obtained by fixing $\beta=0$ and $\gamma$. For the general case, the nominal manifold dimension equals the difference between the number of unknowns and the number of equations, referred to as the \textbf{dimensional deficit}. For equilibria, the dimensional deficit equals $n_p$, which may be reduced to $1$, say, by introducing an additional $n_p-1$ constraints on $x$ and $p$, e.g., that all but one component of $p$ are held fixed.

Under certain \textbf{genericity assumptions}, saddle-node and Hopf bifurcations may be detected and located along one-parameter families of equilibria by monitoring the eigenvalues of the Jacobian $\partial_xf(x,p)$ for crossings of the imaginary axis and then converging on the particular values of $x$ and $p$ associated with the bifurcation. For example, at a saddle-node bifurcation $(x^*,p^*)$, there exists a vector $v^*$ such that $\partial_xf(x^*,p^*)v^*=0$ and $v^{*\mathsf{T}}v^*-1=0$. To continue along a family of saddle-node bifurcations, solutions are obtained from the defining algebraic problem
\[
\mathbb{R}^{2n_x+n_p}\ni\begin{pmatrix}
    x\\p\\v
\end{pmatrix}\mapsto 0=\begin{pmatrix}
    f(x,p)\\\partial_xf(x,p)v\\v^\mathsf{T}v-1
\end{pmatrix}\in\mathbb{R}^{2n_x+1},
\]
which again requires no adaptive updates, since neither the form of the equations, the number of unknowns, nor their meaning change during the analysis. Initialization of this problem is conveniently provided by $(x^*,p^*,v^*)$. In this case, the dimensional deficit equals $n_p-1$ and a one-dimensional solution manifold is obtained by introducing an additional $n_p-2$ constraints on $x$, $p$, and $v$, e.g., that all but two components of $p$ are held fixed.

In the special case that a branch of equilibria $(x(s),p(s))$ is obtained while keeping $n_p-1$ components of $p$ fixed, saddle-node bifurcations coincide with points $(x(s^*),p(s^*))$ with $p'(s^*)=0$, such that $x'(s^*)/\|x'(s^*)\|$ is a candidate for $v^*$. Only in this case may $x'(s^*)$ be used to initialize continuation along a family of saddle-node bifurcations.

At a Hopf bifurcation $(x^*,p^*)$, there exists vectors $v^*$, $w^*$, and $\xi^*$ and a scalar $\omega^*\ne0$, such that $v^{*\mathsf{T}}v^*=1$, $\xi^{*\mathsf{T}}v^*=0$, and $\partial_xf(x^*,p^*)(v^*+\mathrm{i}w^*)=\mathrm{i}\omega^*(v^*+\mathrm{i}w^*)$ or, equivalently, $\partial_xf(x^*,p^*)v^*=-\omega^* w^*$ and $\partial_xf(x^*,p^*)w^*=\omega^* v^*$. To continue along a family of Hopf bifurcations, solutions may be obtained, for example, from the defining algebraic problem
\begin{align*}
\mathbb{R}^{5n_x+n_p+1}\ni\begin{pmatrix}
    x\\p\\k\\v_1\\w_1\\v_2\\w_2
\end{pmatrix}\mapsto 0=\begin{pmatrix}
    f(x,p)\\\partial_xf(x,p)v_1-w_1\\\partial_xf(x,p)v_2-w_2\\
w_1-v_2\\kv_1+w_2\\v_1^\mathsf{T}v_1-1\\\xi^\mathsf{T}v_1
\end{pmatrix}\in\mathbb{R}^{5n_x+2}
\end{align*}
for some vector $\xi$ that picks a unique solution from an infinite family of possible vectors $v_1$. This defining problem requires adaptive updates of $\xi$ during the analysis to ensure that the corresponding Jacobian is of full rank throughout. Initialization when $\xi=\xi^*$ is conveniently provided by $(x^*,p^*,\omega^{*2},v^*,-\omega^*w^*,-\omega^*w^*,-\omega^{*2}v^*)$. In this case, the dimensional deficit equals $n_p-1$ and a one-dimensional solution manifold is obtained by introducing an additional $n_p-2$ constraints on $x$, $p$, and $v$, e.g., that all but two components of $p$ are held fixed.

As advertised in the previous section, the defining problem for Hopf bifurcations shown above may produce solutions that are not Hopf bifurcations. Indeed, by eliminating $v_2$, $w_1$, and $w_2$, we obtain equivalently
\begin{align}
\label{eq:hopf:equiv}
    0=kv_1+J^2v_1,\,0=v_1^\mathsf{T} v_1-1,\,0=\xi^\mathsf{T} v_1, 
\end{align}
where $J=\partial_xf(x,p)$. (Note that the numerical implementation does not perform such elimination since if $J$ has large entries, entries in $J^2$ may be even larger.)  Suppose that a vector $\hat{v}_1$ satisfies the first equation in \eqref{eq:hopf:equiv}. Then $-k$ must be an eigenvalue of $J^2$ with eigenvector $\hat{v}_1$. If $k>0$, then $\hat{v}_1$ cannot be an eigenvector also of $J$ (since $-k<0$). It follows that $\hat{v}_1$ and $J\hat{v}_1$ span a two-dimensional invariant subspace of $J$ for the purely imaginary eigenvalue $\sqrt{-k}$ with corresponding complex eigenvector $\sqrt{-k}\hat{v}_1+J\hat{v}_1$. In particular, there exist locally unique scalars $\alpha$ and $\beta$ such that the vector $\alpha\hat{v}_1+\beta J\hat{v}_1$ satisfies all three of the equations in \eqref{eq:hopf:equiv}.

Suppose, instead, that $k<0$. Then, if $\hat{v}_1$ is not an eigenvector of $J$, $\hat{v}_1$ and $J\hat{v}_1$ span an invariant subspace of $J^2$, which is also an invariant subspace for $J$. It follows that $J$ must have the real eigenvalues $\sqrt{-k}$ and $-\sqrt{-k}$, i.e., the corresponding equilibrium is a neutral saddle. There again exist locally unique scalars $\alpha$ and $\beta$ such that the vector $\alpha \hat{v}_1+ \beta J \hat{v}_1$ satisfies all three of the equations in \eqref{eq:hopf:equiv}.

If, instead, $\hat{v}_1$ is an eigenvector of $J$ (by necessity corresponding to an eigenvalue $\sqrt{-k}$ or $-\sqrt{-k}$), then $\hat{v}_1$ and $J\hat{v}_1$ do not span a two-dimensional invariant subspace for $J$ and there are no constraints on the other eigenvalues of $J$. Orthogonality with respect to $\xi$ here plays the role of the eigenvalue conditions in the preceding two paragraphs. In the CSTR model example, this case corresponds to the solution family \eqref{cstr:non-neutral}. Such a solution family may intersect a family of neutral saddles at a unique branch point, as observed in the preceding section for the CSTR model example.

\subsection{Boundary-value problems}
The defining system for the computation of periodic orbits in Section~\ref{sec:cstr:per} is based on a discretization of the time history using piecewise polynomial collocation \cite{ascher1981collocation}. The details of the implementation for \textsc{coco} are developed in \cite{DS13}, but similar methods are used in \texttt{AUTO} (see \cite{D07} for a detailed description of the discretization) and \texttt{MatCont} \cite{DGK03}. For an autonomous differential equation $x'(t)=f(x(t),p)$ on an arbitrary time interval $[0,T]$, one first rescales time to obtain
\begin{align}
  \label{gen:bvp}
  x'(\tau)=Tf(x(\tau),p)\mbox{ for }\tau\in[0,1].
\end{align}
The function $\tau\mapsto x(\tau)$ is next approximated by a continuous function $\tau\mapsto x_\mathrm{d}(\tau)$ with $L$ polynomial pieces on subintervals defined on a mesh $0=\tau_0<\tau_1<\ldots \tau_L=1$, each of degree $n_\mathrm{deg}$. Finally, let $\tau_{\mathrm{coll},j}$ for $j=1,\ldots,n_\mathrm{deg}$ be \textbf{collocation parameters} defining the discretization method (e.g., Gauss-Legendre nodes \cite{IK94} scaled to a base interval $[0,1]$) and $P_\mathrm{d}$ be the projection that maps the function $\tau\mapsto y(\tau)$ onto the unique (possibly discontinuous) piecewise-polynomial function $\tau\mapsto [P_\mathrm{d}y](\tau)$ that is of degree $n_\mathrm{deg}-1$ on each mesh interval and satisfies
\begin{align*}
  [P_\mathrm{d}y](\tau_{i,j})&=y(\tau_{i,j})\mbox{ for }i=1,\ldots,L,\,
  j=1,\ldots,n_\mathrm{deg}\mbox{, and }\tau_{i,j}=\tau_{i-1}+(\tau_i-\tau_{i-1})\tau_{\mathrm{coll},j}.
\end{align*}
Then, the defining problem for a discretized periodic orbit is given by
\begin{align}
  \label{gen:coll:eq}
  &\mathbb{R}^{n_x(Ln_\mathrm{deg}+1)+n_p+1}\ni\begin{pmatrix}
      x_\mathrm{d}\\p\\T
  \end{pmatrix}\mapsto 0=\begin{pmatrix}
      x_\mathrm{d}'-P_\mathrm{d}\left[Tf(x_\mathrm{d},p)\right]\\x_\mathrm{d}(0)-x_\mathrm{d}(1)\\\langle \xi,x_\mathrm{d}\rangle
  \end{pmatrix}\in\mathbb{R}^{n_x(Ln_\mathrm{deg}+1)+1}
\end{align}
for some function $\xi$ that picks a locally unique solution $x_\mathrm{d}$ from the family of time-shifted solutions $\theta\mapsto x_\mathrm{d}((\cdot)+\theta)$ \cite{DS13,D07}. Initialization requires an initial solution guess $(x_{\mathrm{d},0},p_0,T_0)$ and may also include information about a preferred initial tangent vector $(x'_\mathrm{d}(s),p'(s),T'(s))$ along a family $(x_\mathrm{d}(s),p(s),T(s))$ of such solutions. The dimensional deficit equals $n_p$ and a one-dimensional solution manifold is obtained by introducing an additional $n_p-1$ constraints on $x_\mathrm{d}$, $p$, and $T$, e.g., that all but one component of $p$ are held fixed or that all but two components of $p$ and the period $T$ are held fixed. Each of these choices were made in the analysis of the CSTR model example.

\begin{figure}[t]
  \centering
  \includegraphics[width=0.6\textwidth]{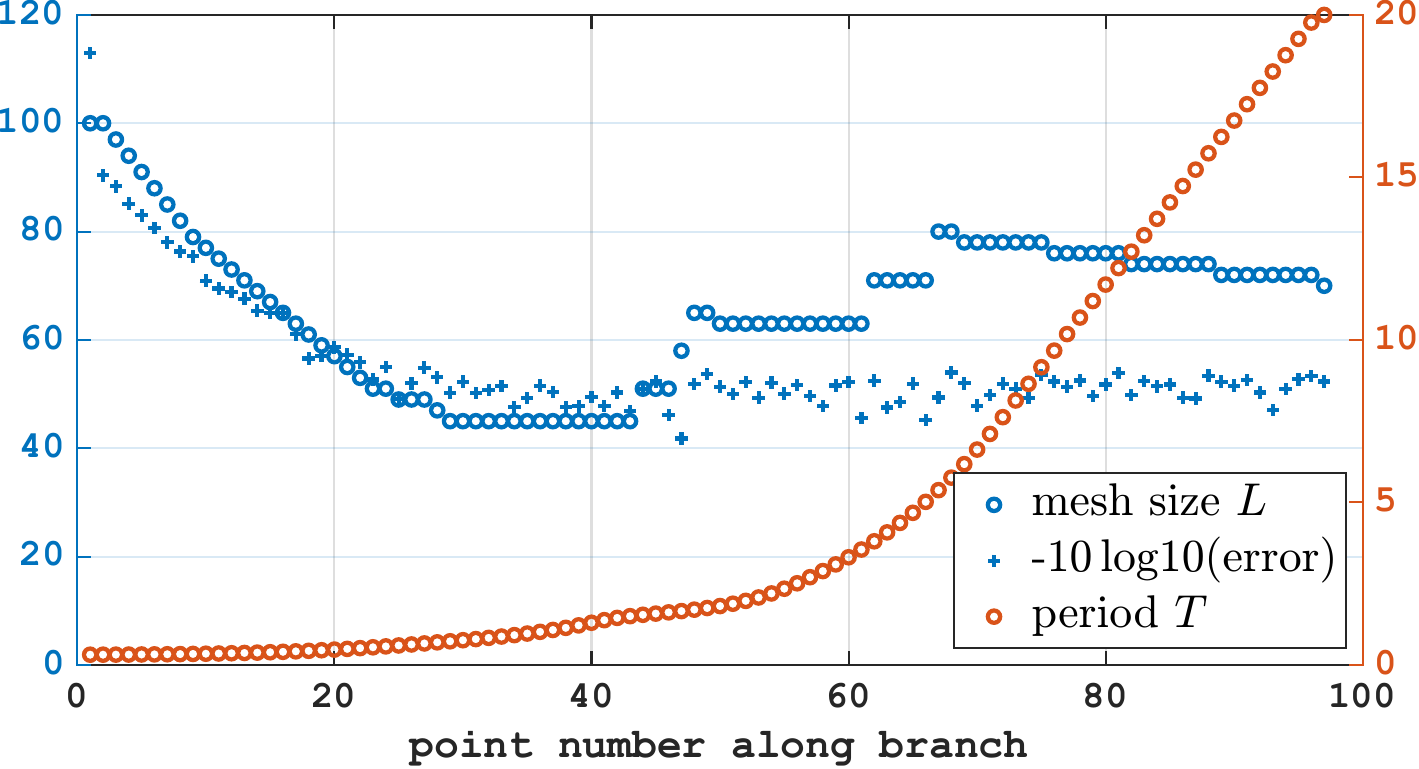}
  \caption[Variable number of unknowns along periodic orbit branches]{Variable number of unknowns along periodic orbit branches shown in Fig.~\ref{fig:homsnic}. Here, the number
$L$ of mesh intervals (called \mcode{NTST} in \textsc{coco}) is automatically adjusted along the curve to ensure that the discretization error does not exceed a given value, keeping the error below $10^{-4}$ (above the horizontal grid line for
$40$ of the left $y$-axis).}
  \label{gen:varpo-illu}
\end{figure}

The defining problem for a periodic orbit requires adaptive updates of $\xi$ to ensure that the problem remains regular. An additional form of adaptive update is obtained if one permits changes to $L$ and the mesh $\{\tau_i\}_{i=1}^{L-1}$ during the analysis, for example, in order to ensure that the discretization error does not exceed a given threshold. As an example, for the initial continuation run along a family of periodic orbits for the model CSTR problem in Section~\ref{sec:cstr:num}, Fig.~\ref{gen:varpo-illu} shows variations in $L$ (as well as in the period and estimated discretization error) with location along the family.

As already remarked upon in Section~\ref{sec:cstr:num}, for the special case of periodic orbits emanating from a generic Hopf bifurcation at $(x^*,p^*)$ and with corresponding vectors $v^*$ and $w^*$ and scalar $\omega^*$, an initial solution guess $x_{\mathrm{d},0}$ may be constructed from a discretization of $t\rightarrow x^*+v^*\cos\omega^* t-w^*\sin\omega^* t$ such that $T_0=2\pi/\omega^*$.

Under certain genericity assumptions, bifurcations may be detected and located along one-parameter families of periodic orbits by monitoring for crossings of the unit circle of the eigenvalues of the monodromy matrix $X(1)$ obtained from the solution to the first variational initial-value problem
\[
\dot{X}=\partial_xf(x,p)X,\,X(0)=I_{n_x}
\]
or, after rescaling and discretization, the solution to the algebraic problem
\[
X_\mathrm{d}'=P_\mathrm{d}[T\partial_xf(x_\mathrm{d},p)X_\mathrm{d}],\,X_\mathrm{d}(0)=I_{n_x},
\]
where $X_\mathrm{d}$ is a continuous, piecewise-polynomial approximant on the identical mesh as $x_\mathrm{d}$. For example, at a saddle-node bifurcation $(x^*_\mathrm{d},p^*,T^*)$ with corresponding monodromy matrix $X_\mathrm{d}^*(1)$, there exists a vector $v^*$ that is orthogonal to $f(x^*_\mathrm{d}(0),p^*)$ and a scalar $b^*$ such that $X_\mathrm{d}^*(1)v^*-v^*=b^*f(x^*_\mathrm{d}(0),p^*)$ and $v^{*\mathsf{T}}v^*-1=0$. To continue along a family of saddle-node bifurcations, solutions are obtained from the defining problem
\begin{align}\label{bvp:sn:defining:problem}
\mathbb{R}^{n_x(2Ln_\mathrm{deg}+3)+n_p+2}\ni\begin{pmatrix}
      x_\mathrm{d}\\p\\T\\y_\mathrm{d}\\b\\v
  \end{pmatrix}\mapsto0=&\,\begin{pmatrix}
      x_\mathrm{d}'-P_\mathrm{d}\left[Tf(x_\mathrm{d},p)\right]\\
      x_\mathrm{d}(0)-x_\mathrm{d}(1)\\
      \langle \xi,x_\mathrm{d}\rangle\\
      y_\mathrm{d}'-P_\mathrm{d}[T\partial_xf(x_\mathrm{d},p)y_\mathrm{d}]\\
      f(x_\mathrm{d}(0),p)^\mathsf{T}v\\
      y_\mathrm{d}(0)-v\\y_\mathrm{d}(1)-v-bf(x_\mathrm{d}(0),p)\\
      v^\mathsf{T}v-1
  \end{pmatrix}\in\mathbb{R}^{n_x(2Ln_\mathrm{deg}+3)+3}.
\end{align}
This problem requires the same adaptive updates of $\xi$, $L$ and $\{\tau_i\}_{i=1}^{L-1}$ as for the defining problem of a periodic orbit. Initialization of this problem is conveniently provided by $(x_\mathrm{d}^*,p^*,T^*,X_\mathrm{d}^*v^*,b^*,v^*)$. In this case, the dimensional deficit equals $n_p-1$ and a one-dimensional solution manifold is obtained by introducing an additional $n_p-2$ constraints on $x$, $p$, and $v$, e.g., that all but two components of $p$ are held fixed.

The defining problem \eqref{bvp:sn:defining:problem} for a saddle-node bifurcation of periodic orbits is guided by the desire to ensure that the monodromy matrix associated with the variational problem has a generalized eigenvector $y_\mathrm{d}(0)$ corresponding to the Floquet multiplier $1$, orthogonal to $\dot x_\mathrm{d}(0)=f(x_\mathrm{d}(0),p)$. In a similar fashion, other generic local bifurcations of periodic orbits, such as \textbf{period-doubling bifurcations} or \textbf{torus bifurcations} have defining problems that ensure that the monodromy matrix has an eigenvector for a desired Floquet multiplier. For period-doubling bifurcations, this is achieved by omitting the fifth component in \eqref{bvp:sn:defining:problem} and replacing the next-to-last component with $y_\mathrm{d}(1)+v$, such that $y_\mathrm{d}(0)$ is an eigenvector corresponding to a Floquet multiplier $-1$. 

For the torus bifurcation, a defining problem is obtained by (i) replacing the unknown function $y_\mathrm{d}$ with two distinct unknown functions $y_\mathrm{d,1}$ and $y_\mathrm{d,2}$, (ii) replacing the unknown vector $v$ with two unknown vectors $v_1$ and $v_2$, (iii) replacing the unknown scalar $b$ with two unknown scalars $a$ and $b$, (iv) replacing the fourth component in \eqref{bvp:sn:defining:problem} with two identical copies for $y_\mathrm{d,1}$ and $y_\mathrm{d,2}$, respectively, and (v) replacing the bottom four components of \eqref{bvp:sn:defining:problem} with the vector
    \[
    \begin{pmatrix}
        y_\mathrm{d,1}(0)-v_1\\
        y_\mathrm{d,2}(0)-v_2\\
        y_\mathrm{d,1}(1)-av_1+bv_2\\
        y_\mathrm{d,2}(1)-a v_2-b v_1\\
        v_1^\tran v_1+v_2^\tran v_2-1\\
        v_1^\tran v_2\\
        a^2+b^2-1
    \end{pmatrix}
    \]
so that $y_\mathrm{d,1}(0)\pm\mathrm{i}y_\mathrm{d,2}(0)$ are eigenvectors of the monodromy matrix corresponding to the Floquet multipliers $a\pm\mathrm{i}b$ on the unit circle.

\subsection{Periodic orbits near homoclinics}
\label{sec:homoclinic}
Although homoclinic connecting orbits may be analyzed in terms of a dedicated defining system \cite{beyn1990numerical}, they may also be analyzed within some desired degree of accuracy by considering the defining system of periodic orbits with fixed and sufficiently large period (while fixing no more than $n_p-2$ components of $p$). This approximation relies on the observation that for a homoclinic orbit and each sufficiently large $T$, there exists a unique periodic orbit of period $T$ close to the homoclinic orbit (up to phase shift). This observation was relied upon in Section~\ref{sec:cstr} where an increasing period along a family of periodic orbits accompanying an accumulation onto a limiting value of a system parameter suggested the existence of a nearby homoclinic orbit.

Two kinds of homoclinic orbits were considered in Section~\ref{sec:cstr}, namely, those associated with a saddle equilibrium and those associated with a saddle-node bifurcation. As pointed out in
\cite{sandstede1997convergence}, under one-parameter variations along families of periodic orbits approaching a homoclinic orbit, the rates at which the period $T$ of nearby periodic orbits grows with the parameter distance to the limiting parameter value differ qualitatively between these two scenarios.

\begin{figure}[ht]
  \centering
  \includegraphics[width=0.75\textwidth]{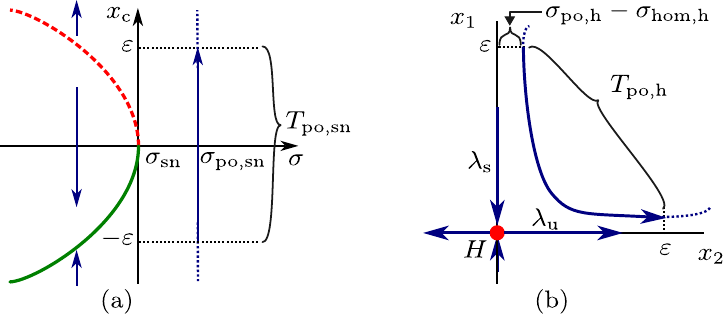}
  \caption[Illustration for asymptotics of period near homoclinic orbits]{Illustration for asymptotics of period near homoclinic orbits: (a) Passage of periodic orbit through $\delta$-neighborhood of saddle-node in extended center manifold.  (b) Passage through $\delta$-neighborhood of hyperbolic saddle $H$ with eigenvalues $\lambda_\mathrm{s}<0<\lambda_\mathrm{u}$. }
  \label{gen:homoclinic-illu}
\end{figure}

Consider, first, the scenario shown in Fig.~\ref{fig:homsnic}, in which a family of periodic orbits accumulate on a connecting orbit to a saddle-node bifurcation equilibrium at a \textbf{saddle-node-on-invariant curve bifurcation} as $\sigma\downarrow\sigma_\mathrm{snic}$. Here, increases in the period of nearby periodic orbits are dominated by a slowing down of the dynamics in the immediate vicinity of the ghost of the saddle-node bifurcation equilibrium as represented schematically in Fig.~\ref{gen:homoclinic-illu}(a). This depicts a size-$\varepsilon$ neighborhood of the saddle-node bifurcation equilibrium in the center manifold (horizontal axis is $\sigma$, vertical axis is the other center manifold direction, called $x_\mathrm{c}$). The parabolic curve in
Fig.~\ref{gen:homoclinic-illu}(a) is the family of equilibria associated with the saddle-node bifurcation (red dashed unstable, and green solid stable). Expanding the vector field to second order about the saddle-node bifurcation,  one concludes that the time for $x_\mathrm{c}$ to increase from $-\varepsilon$ to $\varepsilon$ for $\sigma>\sigma_\mathrm{snic}$ varies as 
\begin{align*}
  \Delta T\sim\frac{1}{\sqrt{\sigma-\sigma_\mathrm{snic}}}
\end{align*}
for $\sigma$ close to $\sigma_\mathrm{snic}$.

In contrast, consider the scenario in which a family of periodic orbits accumulate on a connecting orbit to a saddle equilibrium at a \textbf{homoclinic bifurcation} as $\sigma\downarrow\sigma_\mathrm{hom}$. Here, increases in the period of nearby periodic orbits are dominated by a slowing down of the dynamics in the immediate vicinity of the saddle equilibrium as represented schematically in Fig.~\ref{gen:homoclinic-illu}(b). This depicts a size-$\varepsilon$ neighborhood of the saddle equilibrium, including its stable and unstable manifolds. The distance between the stable manifold and the periodic orbit is proportional to the $\sigma-\sigma_\mathrm{hom}$. If the saddle equilibrium has dominant eigenvalues $\lambda_\mathrm{s}<0<\lambda_\mathrm{u}$ and the connecting orbit is tangential to their corresponding eigenspaces, then the time $\Delta T$ it takes to pass through the neighborhood is determined by
\begin{align*}
  \dot{x}_2=\lambda_\mathrm{u}x_2+O(x_2^2),\,
  x_2(0)=\sigma-\sigma_\mathrm{hom},\,x_2(\Delta T)=\varepsilon,
\end{align*}
such that
\begin{align*}
  \Delta T\sim\frac{-\log(\sigma-\sigma_\mathrm{hom})}{\lambda_\mathrm{u}}
\end{align*}
for $\sigma$ close to $\sigma_\mathrm{hom}$. All things being equal, a large-$T$ periodic orbit is therefore a much more accurate approximation of a nearby homoclinic orbit that connects a saddle equilibrium to itself than one connecting a saddle-node bifurcation equilibrium to itself (see \cite{sandstede1997convergence}).

\section{Genericity and non-genericity}

Several references to genericity in the preceding section are meant to alert the reader to hidden assumptions that underlie the choice of computational algorithm and the interpretation of results. The realization that genericity may be violated in certain scenarios motivates one to consider appropriate modifications to the computational procedure and the way in which bifurcation theory is applied.

Consider, for example, the identification of generic Hopf bifurcations under variations in a single system parameter as points along branches of equilibria from which emanate locally perpendicular (in amplitude/parameter space) branches of periodic orbits according to classical bifurcation theory \cite{K04}. The crossing of the imaginary axis by a pair of complex conjugate eigenvalues of the Jacobian $\partial_xf(x,p)$ is, of course, a necessary condition for the occurrence of such bifurcations. Conclusive identification, however, further requires that the crossing be transversal, that no other eigenvalues equal integer multiples of the triggering pair, and that the first Lyapunov coefficient $\ell_1$ be non-zero. As these conditions hold generically, and assuming no reasons exist to doubt the genericity of the problem at hand, it typically suffices to first detect a crossing, then locate this within desired accuracy, and finally to identify this with a Hopf bifurcation, without additional characterization. This is the approach taken in Section~\ref{sec:cstr}, where $n_x=2$ relieves one from monitoring any additional eigenvalues, and where monitoring of the sign of $\ell_1$ enables delineation between subcritical and supercritical bifurcations. The appropriate identification of the point associated with vanishing $\ell_1$, however, requires computation of the second Lyapunov coefficient $\ell_2$, which was not undertaken in the analysis of the CSTR. In a similar sense, the identification of generic saddle-node bifurcations under variations in a single system parameter as quadratic fold points along branches of equilibria necessitates not only a sign change in a single real eigenvalue of $\partial_xf(x,p)$, as relied upon in the analysis in Section~\ref{sec:cstr}, but also that the Jacobian $(\partial_xf(x,p),\partial_pf(x,p))$ have full rank, something which was not checked in the analysis of the CSTR.

Dynamical systems with \textbf{symmetries} often violate the genericity assumptions relied upon for identifying saddle-node and Hopf bifurcations only by the crossing of the imaginary axis of individual or single pairs of eigenvalues of $\partial_xf(x,p)$. Such violations necessitate a bifurcation theory specific to the nature of the symmetry. This modifies the meaning of genericity for the purpose of identification and requires the formulation of new defining systems for the purpose of continuation.


\subsection{A Network of Chemical Oscillators}
\label{sec:brus}
As an illustration, the reader is invited to consider the eight-dimensional dynamical system
\begin{equation}\label{brus:model}
    \dot{x}_i=A-(B+1)x_i+x_i^2y_i+\epsilon\sum_{j=1}^4(x_j-x_i),\,
    \dot{y}_i=Bx_i-x_i^2y_i+10\epsilon\sum_{j=1}^4(y_j-y_i)
\end{equation}
in terms of the non-negative state $z=(x_1,y_1,\ldots,x_4,y_4)$ and system parameters $A,B,\epsilon\ge 0$, previously introduced in \cite{wulff2006numerical} to demonstrate the approach of the software tool \texttt{sympercon} to exploring systems with discrete symmetry. When $\epsilon=0$, the four identical pairs of equations constitute a \textbf{Brusselator model} of an autocatalytic chemical reaction with $x$ and $y$ describing the concentrations of two intermediate products in the reaction chain. For $\epsilon\ne0$, these pairs of equations are coupled in a fully connected topology. Consequently, for all values of $\epsilon$, the full model is \textbf{equivariant} under arbitrary permutations of indices. In other words, a solution of \eqref{brus:model} can be transformed into another solution by permuting the subscripts $i$ arbitrarily.

A fully symmetric solution is one for which each of the four pairs of states exhibits the identical time history. As an example, a fully symmetric equilibrium occurs when $x_i=A$ and $y_i=B/A$ for all $i$. Its stability is to linear approximation characterized by the eigenvalues of the Jacobian matrix
\[
J:=I_4\otimes \begin{pmatrix*}[r]B-1 & A^2\\-B & -A^2\end{pmatrix*}+\epsilon(\mathbb{1}_4-4I_4)\otimes \begin{pmatrix*}[r]1& 0 \\0 & 10\end{pmatrix*},
\]
where $\mathbb{1}_4$ is the $4\times4$ matrix with all entries equal to $1$.
Here, permutation symmetry of the model equations implies that $(\Pi\otimes I_2) J=J(\Pi\otimes I_2)$ for an arbitrary permutation matrix $\Pi$ and, consequently, that if $(\kappa,v)$ is an eigenpair of $J$, then so is $(\kappa,(\Pi\otimes I_2)v)$.

With $\delta:=B-1-A^2$, the eigenvalues of $J$ evaluate to
\[
\kappa_{1,\pm}=\frac{\delta\pm\sqrt{\delta^2-4A^2}}{2}
\]
with eigenvectors of the form
\[
(k_1\varphi,\varphi,k_1\varphi,\varphi,k_1\varphi,\varphi,k_1\varphi,\varphi)^\mathsf{T}
\]
for arbitrary $\varphi$ and some fixed scalar $k_1$, and
\[
\kappa_{2,\pm}=\frac{\delta-44\epsilon\pm\sqrt{(\delta+36\epsilon)^2-4A^2(1-36\epsilon)}}{2}
\]
with eigenvectors of the form
\[
(k_2(\varphi_1+\varphi_2+\varphi_3),\varphi_1+\varphi_2+\varphi_3,-k_2\varphi_1,-\varphi_1,-k_2\varphi_2,-\varphi_2,-k_2\varphi_3,-\varphi_3)^\mathsf{T}
\]
for arbitrary $\varphi_1,\varphi_2,\varphi_3$ and some fixed scalar $k_2$. While unexpected for a generic system, equivariance here implies a high-dimensional eigenspace for $\kappa_{2,\pm}$.

For fixed $A$ and $\delta$, such that $11A^2(11-9\delta)>100\delta^2$, and under variations in $\epsilon$, the second pair of eigenvalues, $\kappa_{2,\pm}$, crosses the imaginary axis transversally for $\epsilon=\epsilon^*:=\delta/44$ at $\pm \mathrm{i}\omega$ with
\[
\omega=\frac{\sqrt{11A^2(11-9\delta)-100\delta^2}}{11}\mbox{ and }k_2=\frac{11A^2}{11i\omega-11A^2-10\delta}.
\]
The generic Hopf bifurcation theorem required as one of its genericity conditions that the matrix $J|_{\epsilon=\epsilon^*}-\i\omega I_8$ have a one-dimensional complex nullspace. This condition fails here as the permutation symmetry of \eqref{brus:model} leads to a three-dimensional complex nullspace of $J|_{\epsilon=\epsilon^*}-\mathrm{i}\omega I_8$.  

The \textbf{equivariant Hopf bifurcation theorem} permits construction of families of periodic orbits that emanate from such a bifurcation point (see \cite{golubitsky1985hopf,golubitsky2012singularities} for the general theory in the presence of symmetry, \cite{dias2009hopf} for systems with general permutation symmetry, and \cite{swift1988hopf} for a study of the case of four oscillators coupled in a ring topology). Periodic orbits in a system with symmetry are classified according to their \textbf{spatiotemporal symmetry}, here represented by a permutation composed of $K$ subcycles and an integer  \textbf{cycle shift} for each subcycle. For example, the notation $(1\,2\,4\,3)_4$ applies to a periodic orbit of period $T$ that satisfies the conditions
\begin{equation}
    \label{brus:(1,2,4,3)_4}
    \begin{pmatrix}x_1\\y_1\end{pmatrix}(t)=\begin{pmatrix}x_2\\y_2\end{pmatrix}\left(t+\frac{T}{4}\right)=\begin{pmatrix}x_4\\y_4\end{pmatrix}\left(t+\frac{T}{2}\right)=\begin{pmatrix}x_3\\y_3\end{pmatrix}\left(t+\frac{3T}{4}\right),
\end{equation}
with $K=1$ and single cycle shift $4$, while $(1\,2)_2(3\,4)_1$ applies to a periodic orbit of period $T$ that satisfies the conditions
\begin{equation}
    \label{brus:(1,2)_2(3,4)_1}
    \begin{pmatrix}x_1\\y_1\end{pmatrix}(t)=\begin{pmatrix}x_2\\y_2\end{pmatrix}\left(t+\frac{T}{2}\right),\,\begin{pmatrix}x_3\\y_3\end{pmatrix}\left(t\right)=\begin{pmatrix}x_4\\y_4\end{pmatrix}\left(t+\frac{T}{1}\right)=\begin{pmatrix}x_4\\y_4\end{pmatrix}\left(t\right)
\end{equation}
with $K=2$ and cycle shifts $2$ and $1$. Let $\ell_k$ equal the cycle shift of the subcycle containing index $k$, such that $\ell=(\ell_1,\ell_2,\ell_3,\ell_4)$ is the corresponding \textbf{shift vector}.
Then, $\ell$ is $(4,4,4,4)$ for the spatiotemporal symmetry $(1\,2\,4\,3)_4$ in \eqref{brus:(1,2,4,3)_4}, while $\ell=(2,2,1,1)$ for the spatiotemporal symmetry $(1\,2)_2(3\,4)_1$ in \eqref{brus:(1,2)_2(3,4)_1}. When a subcycle contains only a single index and the corresponding cycle shift equals $1$, we sometimes omit the subcycle from the representation.

For the existence of a branch of periodic orbits of some spatiotemporal symmetry with permutation matrix $\Pi$ and shift vector $\ell$, the {equivariant Hopf bifurcation theorem then requires checking if the complex dimension of the nullspace of the augmented matrix $J_{\Pi,\ell}$ below equals $1$, i.e., if
\begin{align}\label{brus:equiv:hopf}
\dim_\mathbb{C}\ker J_{\Pi,\ell}&=1\mbox{,\quad where\quad}J_{\Pi,\ell}=\begin{bmatrix}
    J|_{\epsilon=\epsilon^*}-\i\omega I_8\\
    \Pi\otimes I_2-E_\ell\otimes I_2
\end{bmatrix}
\end{align}
and $E_\ell$ denotes the diagonal $4\times 4$ matrix with $E_{\ell,(k,k)}=\e^{2\pi\i/\ell_k}$ in the diagonal entries. For permutations that satisfy this condition, there follows the existence of a branch of periodic orbits $z(t)$ of period $T$, which (i) emanates from the bifurcation point $(z^*,\epsilon^*)$ under either increasing or decreasing variations in $\epsilon$, (ii) satisfies the spatiotemporal symmetry
\begin{equation}
\label{brus:symmx}
\left[(\Pi\otimes I_2)z\right]_k(t)=z_k(t+T/\ell_k),
\end{equation}
and, for $\epsilon$ close to $\epsilon^*$, (iii) takes the form $z(t)-z^*\sim v\cos\omega t-w\sin\omega t$, where $v+\i w$ spans the nullspace of $J_{\Pi,\ell}$, satisfying
\begin{equation}
\label{brus:eigconds}
    J|_{\epsilon=\epsilon^*}v=-\omega w,\,J|_{\epsilon=\epsilon^*}w=\omega v\mbox{, and } (\Pi\otimes I_2)(v+\mathrm{i}w)=E_\ell(v+\mathrm{i}w).
\end{equation}
Indeed, in this limit $T\approx 2\pi/\omega$ and
\begin{equation}
\label{brus:symmv}
\left[(\Pi\otimes I_2)(v\cos\omega t-w\sin\omega t)\right]_k=v_k\cos(\omega t+2\pi/\ell_k)-w_k\sin(\omega t+2\pi/\ell_k).
\end{equation}

As an example, let $\Pi$ be a permutation matrix such that $\Pi^4=I_4$ and $\Pi^k\ne I_4$ for $k<4$, i.e., such that $K=1$ and $\ell=(4,4,4,4)$. If $p_j$ denotes the number of applications of $\Pi$ required to map index $j$ to $1$, $\Pi\otimes I_2$ must have an eigenvalue $\i$ with eigenvectors of the form $(\mu, \nu, \i^{p_2}\mu,\i^{p_2}\nu,\i^{p_3}\mu,\i^{p_3}\nu,\i^{p_4}\mu,\i^{p_4}\nu)^\mathsf{T}$. Conditions \eqref{brus:eigconds} are then satisfied by $v$ and $w$ given by the real and imaginary parts of a vector of this form if $\mu=rA^2$ and $\nu=r(i\omega-A^2-10\delta/11)$ for any arbitrary scalar $r$.

For a quantitative example, suppose (as in \cite{wulff2006numerical}) that $A=2.0$ and $B=5.9$. In this case, the equivariant Hopf bifurcation occurs for $1000\epsilon^*\approx 20.45$ with $\omega^*\approx 0.62$ and \eqref{brus:equiv:hopf} is satisfied by the spatiotemporal symmetry $(1\,2\,4\,3)_4$. The equivariant Hopf bifurcation theorem then predicts the existence of a family of periodic orbits emanating from the bifurcation point with limiting period $2\pi/\omega^*\approx 10.12$ and satisfying the conditions \eqref{brus:(1,2,4,3)_4}. The bifurcation diagram  near the equivariant Hopf bifurcation in Fig.~\ref{fig:brus:hopf} shows the corresponding branch labeled $P_1$ with a typical time profile in the right panel. 
\begin{figure}[ht]
    \centering
\includegraphics[width=0.95\linewidth]{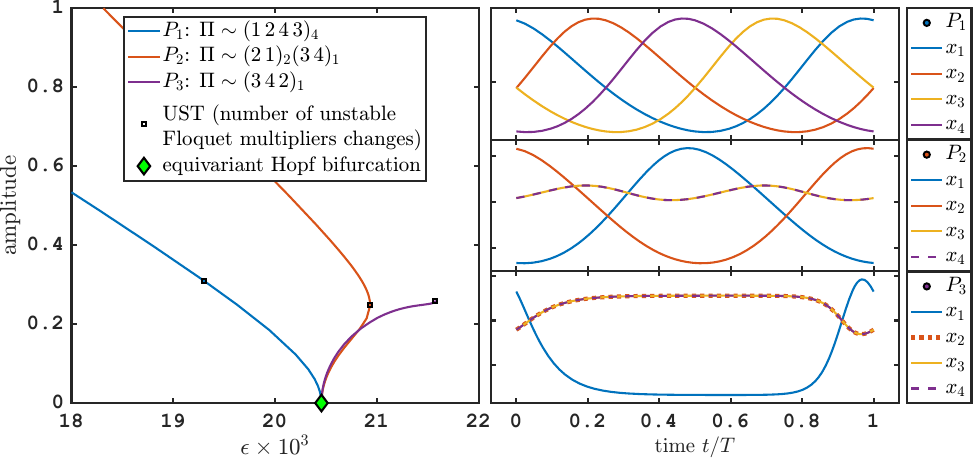}
\caption{Branches of periodic orbits emanating from an equivariant Hopf bifurcation for \eqref{brus:model} with different spatiotemporal symmetries. The subscript $k$ for each cycle indicates the time shift $T/k$ applied for each permutation. The time profiles of $x_i$ in the right panel illustrate the symmetry of the orbits. Notably, the branch denoted by $P_3$ appears to terminate at a homoclinic bifurcation, while the other two branches continue outside of the range shown. Other parameters: $A=2$, $B=5.9$.}
    \label{fig:brus:hopf}
\end{figure}

At the equivariant Hopf bifurcation, condition \eqref{brus:equiv:hopf} is also satisfied by other spatiotemporal symmetries, such as $(1\,2)_2(3\,4)_1$ and $(1)_1(2\,3\,4)_1$ (abbreviated to $(2\,3\,4)_1$). The equivariant Hopf bifurcation theorem consequently predicts the existence of branches of periodic orbits satisfying the conditions \eqref{brus:(1,2)_2(3,4)_1} and
\[
\begin{pmatrix}
    x_2\\y_2
\end{pmatrix}(t)=\begin{pmatrix}
    x_3\\y_3
\end{pmatrix}(t)=\begin{pmatrix}
    x_4\\y_4
\end{pmatrix}(t),
\]
respectively. This is verified by the typical time profiles shown in the right panel of Fig.~\ref{fig:brus:hopf} for the branches of periodic orbits labeled $P_2$ and $P_3$ in the left panel. Each of these three branches also coexists with its images under arbitrary permutations of indices.

In the screen output below, continuation is performed further along the branch $P_1$ of periodic orbits under variations in $\lambda=1000\epsilon$ and with initial solution guess obtained from the preceding theory. The extended branch is shown as $P_1$ in Fig.~\ref{fig:brus:sympo}.
\begin{lstlisting}[language=coco-highlight,frame=lines]
...   LABEL  TYPE        lambda    po.period    amplitude
...       1  EP      2.0455e+01   1.0125e+01   7.9991e-04
...       2  UST     1.9302e+01   8.2014e+00   3.0964e-01
...       3  UST     1.7002e+01   6.4430e+00   6.9931e-01
...       4  BP      1.7002e+01   6.4430e+00   6.9931e-01
...       5  UST     2.6735e+00   5.0410e+00   4.7702e+00
...       6  BP      2.6734e+00   5.0410e+00   4.7702e+00
...       7  EP      1.0000e+00   5.0618e+00   5.6015e+00
\end{lstlisting}
Here, bifurcation points associated with changes in the number of Floquet multipliers outside the unit circle are denoted by the label \mcode{UST}. For the bifurcation at $\lambda\approx 19.30$, two complex conjugate Floquet multipliers cross the unit circle at a generic \textbf{Neimark-Sacker bifurcation}.
\begin{figure}[ht]
  \centering
  \includegraphics[width=0.9\textwidth]{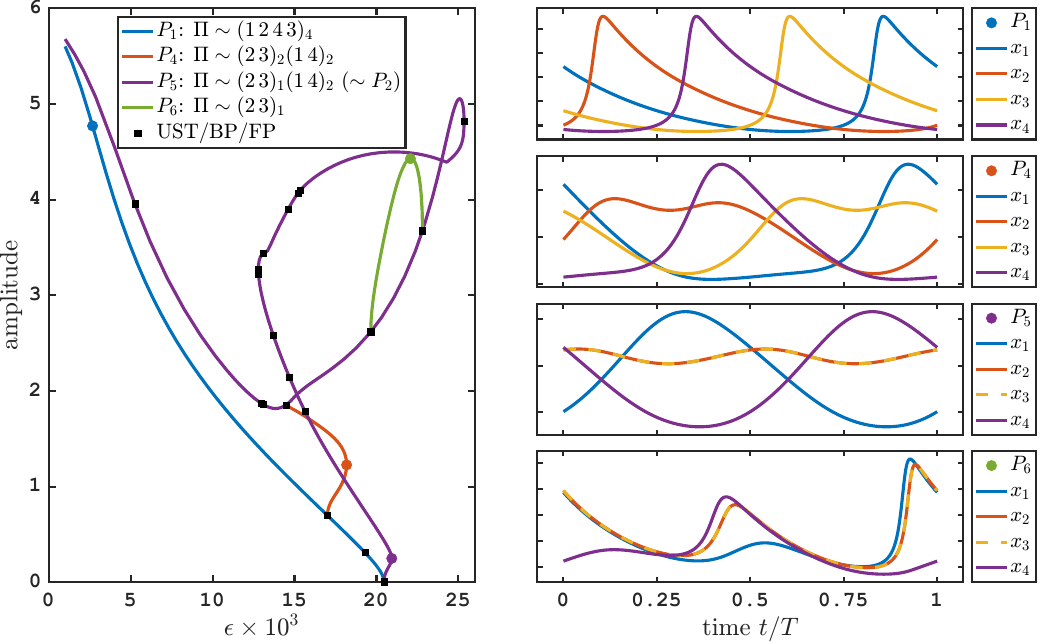}
  \caption[Symmetry breaking of periodic orbits with $(1243)$ symmetry]{Family $P_1$ of periodic orbits of \eqref{brus:model} with $(1\,2\,4\,3)_4$ symmetry, with further symmetry-breaking bifurcations, resulting in families $P_4$, $P_5$, $P_6$ with other symmetries as indicated in the legend and illustrated with time profiles.}
  \label{fig:brus:sympo}
\end{figure}

For each of the bifurcations at $\lambda\approx 17.00$ and $\lambda\approx 2.67$, on the other hand, a single eigenvalue crosses the unit circle at $1$. In contrast to the non-symmetric, generic case, these are not saddle-node bifurcations. Indeed, as the points are also detected as branch points, denoted by \mcode{BP}, and under appropriate genericity assumptions, they are expected to be \textbf{symmetry-breaking flip-pitchfork bifurcations}, an expectation that is borne out by further computational analysis. From each emanate secondary branches of periodic orbits that are symmetric with respect to the spatiotemporal symmetry associated with the permutation $(1\,4)_2(2\,3)_2$, i.e., such that
\[
\begin{pmatrix}x_1\\y_1\end{pmatrix}(t)=\begin{pmatrix}x_4\\y_4\end{pmatrix}\left(t+\frac{T}{2}\right),\,\begin{pmatrix}x_2\\y_2\end{pmatrix}\left(t\right)=\begin{pmatrix}x_3\\y_3\end{pmatrix}\left(t+\frac{T}{2}\right).
\]
The screen output below is obtained from continuation along one such branch, shown as $P_4$ in Fig.~\ref{fig:brus:sympo}. 
\begin{lstlisting}[language=coco-highlight,frame=lines]
...   LABEL  TYPE        lambda    po.period    amplitude
...       1  EP      1.7002e+01   6.4430e+00   6.9931e-01
...       2  UST     1.8178e+01   7.3991e+00   1.2269e+00
...       3  FP      1.8178e+01   7.3991e+00   1.2269e+00
...       4  FP      1.4473e+01   5.9863e+00   1.8534e+00
...       5  BP      1.4473e+01   5.9863e+00   1.8534e+00
...       6  UST     1.8178e+01   7.3991e+00   1.2271e+00
...       7  FP      1.8178e+01   7.3991e+00   1.2271e+00
...       8  FP      1.7002e+01   6.4430e+00   6.9938e-01
...       9  BP      1.7002e+01   6.4430e+00   6.9934e-01
...      10  UST     1.8178e+01   7.3991e+00   1.2272e+00
...      11  FP      1.8178e+01   7.3991e+00   1.2272e+00
...      12  BP      1.4473e+01   5.9863e+00   1.8538e+00
...      13  FP      1.4473e+01   5.9863e+00   1.8538e+00
...      14  FP      1.8178e+01   7.3991e+00   1.2271e+00
...      15  UST     1.8178e+01   7.3991e+00   1.2271e+00
...      16  EP      1.7021e+01   6.4547e+00   7.3790e-01
\end{lstlisting}
Along this branch, at the points denoted by \mcode{UST}, a single Floquet multiplier crosses the unit circle at $1$. These coincide with fold points, denoted by \mcode{FP}, and constitute generic saddle-node bifurcations of periodic orbits. In contrast, at the fold point coincident with a branch point, a single Floquet multiplier touches, but does not cross, the unit circle at $1$. The second of these coincides with the symmetry-breaking bifurcation along the primary branch. In contrast, as is clear from the time profiles shown in Fig.~\ref{fig:brus:sympo}  (see  $P_5$) for the periodic orbit at the first branch point, here the orbit satisfies the additional symmetry $x_2(t)=x_3(t),\,y_2(t)=y_3(t)$. Analogously, the periodic orbit at the third branch point satisfies the additional symmetry $x_1(t)=x_4(t),\,y_1(t)=y_4(t)$. Both points constitute \textbf{symmetry-increasing bifurcations}.

The screen output below shows the result of continuation from the first of these points along a secondary branch of periodic orbits that satisfy the spatiotemporal symmmetry
\[
\begin{pmatrix}x_1\\y_1\end{pmatrix}(t)=\begin{pmatrix}x_4\\y_4\end{pmatrix}\left(t+\frac{T}{2}\right),\,\begin{pmatrix}x_2\\y_2\end{pmatrix}\left(t\right)=\begin{pmatrix}x_3\\y_3\end{pmatrix}\left(t+\frac{T}{2}\right)=\begin{pmatrix}x_3\\y_3\end{pmatrix}\left(t\right).
\]
\begin{lstlisting}[language=coco-highlight,frame=lines]
...   LABEL  TYPE        lambda    po.period    amplitude
...       1  EP      1.4473e+01   5.9863e+00   1.8536e+00
...       2  BP      1.9650e+01   7.2452e+00   2.6173e+00
...       3  UST     1.9650e+01   7.2452e+00   2.6173e+00
...       4  BP      2.2810e+01   8.5430e+00   3.6744e+00
...       5  UST     2.2810e+01   8.5430e+00   3.6743e+00
...       6  UST     2.5346e+01   1.0683e+01   4.8199e+00
...       7  FP      2.5346e+01   1.0683e+01   4.8199e+00
...       8  BP      1.5355e+01   1.2014e+01   4.0979e+00
...       9  UST     1.5355e+01   1.2014e+01   4.0979e+00
...      10  UST     1.5209e+01   1.2009e+01   4.0666e+00
...      11  BP      1.5209e+01   1.2009e+01   4.0666e+00
...      12  BP      1.4616e+01   1.1820e+01   3.9017e+00
...      13  UST     1.4616e+01   1.1820e+01   3.9017e+00
...      14  EP      1.4000e+01   1.1364e+01   3.6991e+00
\end{lstlisting}
(Recall that $\lambda=\epsilon\times10^3$.) Branch switching to $P_6$ in Fig.~\ref{fig:brus:sympo} from the first of the branch points reduces symmetry again. Along this branch $x_2(t)=x_3(t),\,y_2(t)=y_3(t)$ but the spatiotemporal symmetry between $(x_1,y_1)$ and $(x_4,y_4)$ has been lost, as indicated by the cycle notation $(2\,3)_1$ (short for $(2\,3)_1(1)_1(4)_1$) in the legend and the typical time profiles.

\subsection{A cautionary note and its resolution}
As seen in the preceding section, in a system with symmetry, the types of degeneracies encountered along branches of equilibria and periodic orbits are different from what might be expected in generic systems (such as the CSTR example discussed in Section~\ref{sec:cstr}). For example, branches of periodic orbits of different types of symmetry form a network with intersections at symmetry-breaking and symmetry-increasing bifurcations where crossings by Floquet multipliers of the unit circle are coincident with a loss of rank of the Jacobian of the defining system.

Such isolated rank loss makes computation along branches through branch points unreliable, because close to a branch point the solution is not unique. As a result, the analysis may switch uncontrollably near the branch point to any of the intersecting branches. This is exacerbated by the fact that a spatiotemporal symmetry of a periodic orbit such as \eqref{brus:symmx} is only approximately valid for the discretized approximant of the orbit. The reason for the imperfection is that meshes of time points may not be preserved by the time shifts $T/\ell_k$.

Computational tools designed to account for symmetries, such as \texttt{symcon} \cite{gatermann1991symbolic} and \texttt{sympercon} \cite{wulff2006numerical,wulff2010sympercon}, avoid the need to analyze near-singular problems by modifying the defining systems to explicitly account for the desired symmetry constraints. For example, \texttt{sympercon} eliminates the consequences of non-uniqueness and discretization near symmetry-breaking bifurcations by replacing the periodic-orbit boundary-value problem on the interval $[0,T]$ with a boundary-value problem on a shorter interval, say $[0,T/\ell]$, with boundary conditions that impose the desired symmetry. This approach not only avoids singularities associated with symmetry-breaking bifurcations, but also saves computational effort. The required modifications to the defining system, however, depend on the symmetry type.

An alternative to the reduction methodology of \texttt{symcon} and \texttt{sympercon} is to append symmetry constraints to the generic defining systems. This approach is particularly well suited to \textsc{coco}, where defining systems are naturally constructed in stages. For example, to continue branches of periodic orbits satisfying \eqref{brus:symmx} without suffering the consequences of non-uniqueness and discretization near branch points, one can impose one or several components of the symmetry constraints \eqref{brus:symmx} at finitely many instances on $[0,T]$. Although these additional conditions are nominally redundant along the solution branch, they help regularize\footnote{If the defining system with the redundant constraints takes the form $\Phi(u)=0$ such that $\Phi(u^*)=0$ and $\partial_u\Phi(u^*)$ has full column rank, then there exists a matrix $S$ such that the defining system $\Phi(u)+Sw=0$ is regular at $(u,w)=(u^*,0)$. The matrix $S$ can be computed on the fly and updated in each continuation step. This functionality is available in the release of \textsc{coco} accompanying this chapter.} the defining system. While this disables detection of symmetry-breaking bifurcations as branch points, these may still be detected by a Floquet multiplier crossing the unit circle at $1$. 

To illustrate the effect of adding constraints, the screen output below shows the result of continuation along the $P_1$ branch obtained by adding $8$ redundant constraints to the defining system, choosing to impose \eqref{brus:(1,2,4,3)_4} at a selection of times $t$:
\begin{align*}
x_1(0)&=x_2(T/4),&
x_2(0)&=x_4(T/4),&
x_4(0)&=x_3(T/4),&
x_3(0)&=x_1(T/4),\\
y_1(0)&=y_2(T/4),&
y_2(0)&=y_4(T/4),&
y_4(0)&=y_3(T/4),&
y_3(0)&=y_1(T/4).
\end{align*}
There is some arbitrariness in the choice of time points where redundant conditions are imposed and in the number of conditions. The conditions regularize the problem if their nullspace intersects the nullspace of the linearization of the defining system (in this case \eqref{gen:coll:eq}) trivially in symmetry-breaking points.
\begin{lstlisting}[language=coco-highlight,frame=lines]
...   LABEL  TYPE        lambda    po.period    amplitude
...       1  EP      2.0455e+01   1.0125e+01   6.3564e-04
...       2  UST     2.0455e+01   1.0125e+01   7.9895e-04
...       3  UST     2.0455e+01   1.0125e+01   7.9895e-04
...       4  UST     1.9302e+01   8.2014e+00   3.0933e-01
...       5  UST     1.7002e+01   6.4430e+00   6.9930e-01
...       6  EP      1.5191e+01   5.7379e+00   1.0000e+00
\end{lstlisting}
In contrast to the preceding screen output, the point at $\lambda\approx17$ that was previously detected as a branch point (\texttt{BP}) is now only detected by a change in the number of Floquet multipliers outside the unit circle (\texttt{UST}).

Branch points that coincide with symmetry-increasing bifurcations cannot be avoided by either methodology. At best, one may anticipate some singularities and prevent uncontrolled switches to other branches of solutions by monitoring the deviation of the solution from satisfying such additional spatiotemporal symmetries. Once the character of such additional symmetries is identified (by inspection of the solution), continuation may proceed along the secondary branch using either of the approaches described above.

\subsection{Tracking bifurcations}
\begin{figure}
    \centering
    \includegraphics[width=0.95\textwidth]{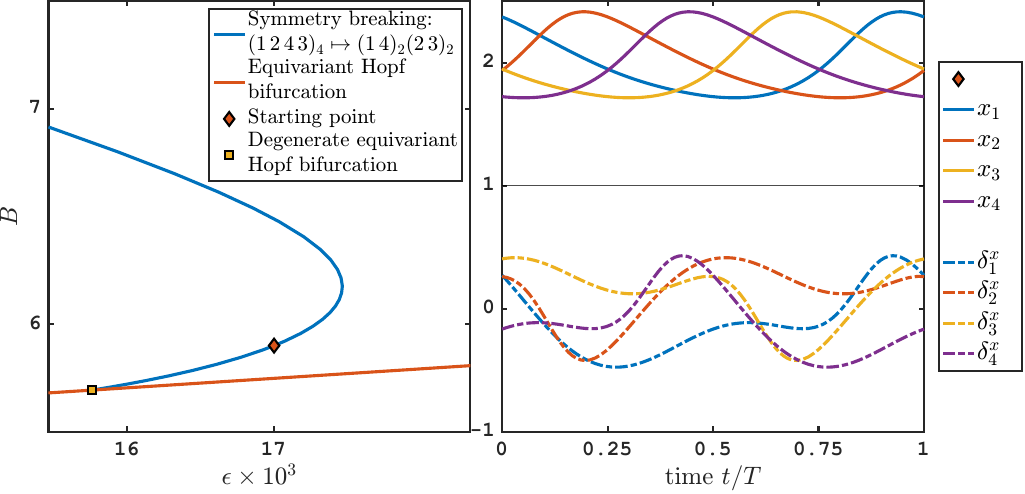}
    \caption{Locus of equivariant Hopf bifurcation and symmetry breaking bifurcation in the $(\epsilon,B)$ two-parameter plane for the model \ref{brus:model}. Time profiles of solution components $x_i$ and variational components $\delta^x_i$ show the symmetries connected by the symmetry-breaking bifurcation.}
    \label{fig:enter-label}
\end{figure}The idea of adding symmetry constraints on top of generic defining systems also applies to the regularized tracking of bifurcations of equilibria and periodic orbits in systems with symmetry. For example, to continue an equivariant Hopf bifurcation associated with a spatiotemporal symmetry $(\Pi,\ell)$, as in \eqref{brus:equiv:hopf}, one may simply append the $n_x$ complex equations (equivalent to $2n_x$ real equations)
\begin{align}
\label{brus:symmhopf}
\begin{aligned}
0&=(\Pi\otimes I_2-E_\ell\otimes I_2)(\sqrt{k}v_1-\mathrm{i}Jv_1),    
\end{aligned}
\end{align}
to the defining system for tracking generic Hopf bifurcation points.

The screen output below shows the result of continuing along a branch of equivariant Hopf bifurcations for \ref{brus:model} from an initial solution guess given by the fully symmetric equilibrum with $(A,B)=(2,5.9)$ under simultaneous variations in $B$ and $\lambda=\epsilon\times10^3$ and for fixed $A$.
\begin{lstlisting}[language=coco-highlight,frame=lines]
...  LABEL  TYPE        lambda            B            k
...      1  EP      2.0455e+01   5.9000e+00   3.8512e-01
...      2  EP      1.5000e+01   5.6600e+00   1.4800e+00
...  LABEL  TYPE        lambda            B            k
...      3  EP      2.0455e+01   5.9000e+00   3.8512e-01
...      4  MX      2.2263e+01   5.9796e+00   1.0422e-03
\end{lstlisting}
We monitor the quantity $k$ (which equals the square of the imaginary part of the critical eigenvalues) along the branch and observe a lack of convergence when $k$ approaches zero near $(\lambda,B)\approx(22.2,5.98)$, due to the non-differentiability in this limit of \eqref{brus:symmhopf}.

For symmetry-breaking bifurcations along branches of periodic orbits with spatiotemporal symmetry \eqref{brus:symmx}, a Floquet multiplier of geometric multiplicity $1$ or higher crosses the unit circle at $1$. 
It follows that continuation of branches of such symmetry-breaking bifurcations is achieved by appending one or several components of \eqref{brus:symmx} at finitely many instances $t\in[0,T]$ to the defining system for tracking generic saddle-node bifurcations of periodic orbits with right-hand side \eqref{bvp:sn:defining:problem}.

The screen output below shows the result of continuing along a branch of symmetry-breaking bifurcations with initial solution guess obtained from the second \mcode{UST} point along the branch of periodic orbits of symmetry \eqref{brus:(1,2,4,3)_4} under simultaneous variations in $\lambda=\epsilon\times10^3$ and $B$ (keeping $A=2$ fixed).
\begin{lstlisting}[language=coco-highlight,frame=lines]
LABEL  TYPE        lambda            B    po.period    amplitude
    1  EP      1.7002e+01   5.9000e+00   6.4431e+00   6.9931e-01
    2  FP      1.5762e+01   5.6935e+00   5.4424e+00   3.3284e-05
    3  EP      1.5780e+01   5.6956e+00   5.4526e+00   7.0158e-02

LABEL  TYPE        lambda            B    po.period    amplitude
    4  EP      1.7002e+01   5.9000e+00   6.4431e+00   6.9931e-01
    5  FP      1.7467e+01   6.1746e+00   7.6670e+00   1.0523e+00
    6  EP      1.5463e+01   6.9152e+00   8.9513e+00   1.9987e+00
\end{lstlisting}
The second point labelled \texttt{FP} is an extremum of $\lambda$ (and $\epsilon$) along the curve, not associated with a bifurcation of higher codimension. 

Depending on the symmetry of the periodic orbit, it may have Floquet multipliers of  geometric multiplicity generically higher than $1$. For example, the periodic orbit family $P_3$ shown in Fig.~\ref{fig:brus:hopf} has double Floquet multipliers, which may cross the unit circle simultaneously. For this case the \textbf{Equivariant Branching Lemma} provides a criterion for finding branches of periodic orbits that emerge from the higher-multiplicity symmetry breaking \cite{golubitsky2012singularities}.

\section{Concluding discussion}
Section~\ref{sec:cstr} showcased bifurcation analysis for the classical use case: the CSTR problem is a low-dimensional ODE depending on multiple parameters, which shows a rich but finite range of qualitatively different phase portraits for different parameters. The parameters for each phase portrait are determined by tracking equilibria or periodic orbits, and their bifurcations. The end result is a map classifying different qualitative behavior depending on parameters, such as Figure~\ref{fig:cstr:bif2d}.

A strength of this approach is that it can be applied to a wide variety of problems using a small number of nonlinear algebraic problems, namely the defining systems for equilibria, periodic orbits and a few \emph{generic codimension}-$1$ bifurcations. Our section~\ref{sec:defining_systems} outlined most of the defining algebraic problems in the form in which they are implemented in \textsc{coco}. Comprehensive overviews of defining algebraic problems are given by \cite{MDO09} and textbooks such as \cite{G00,K04}. The set of defining algebraic problems for classical bifurcation analysis as performed here for the CSTR problem has been generalized to differential equations with delay and renewal equations in \textsc{dde-biftool} \cite{ELR02} (see \cite{ddebiftoolmanual} for manual) and \textsc{knut} \cite{RS07}. These generalizations use the same underlying numerical methods as \textsc{coco}, \textsc{auto} or \textsc{matcont}, such as collocation discretizations for the periodic orbits. However, they require different numerical methods for the computation of stability and eigenvalues, and expressions for variational problems and eigenvalues or eigenvectors look different \cite{breda2005pseudospectral}.

The discretization of periodic orbits with collocation is especially suitable for time profiles with strongly non-harmonic shape such as those shown in Figure~\ref{fig:homsnic}(d) with sharp peaks. Some problem classes, such as vibration analysis for mechanical structures, consider scenarios with simple (e.g., harmonic) forcing but complex spatial shapes. For these problems a general-purpose continuation method for continuation of periodic responses benefits from the spectral convergence of harmonic projections, such that harmonic balance with very few \emph{modes} results in sufficient accuracy \cite{cochelin2009high,detroux2015harmonic}. For example, the vibrations of a satellite structure with $37$ degrees of freedom, used as an illustration in \cite{detroux2015harmonic} could be analyzed using just $5$ harmonics (comparing to $9$ harmonics to confirm accuracy).
The harmonic balance method is also suitable for cases where the right-hand side can only be evaluated with low accuracy but a high sampling rate (e.g., through experimental measurements), because the projection on harmonics has an averaging effect on random disturbances \cite{schilder2015experimental}. For this reason, it has been combined with feedback control to track bifurcations and unstable periodic responses in vibration experiments \cite{barton2017control,renson2016robust}. See \cite{raze2024experimental} for a review of methods for bifurcation analysis in controlled experiments.

For high-dimensional problems arising from discretizations of PDEs in fluid dynamics time-stepper based methods, combined with matrix-free linear algebra for evaluation of Jacobian matrices and linear equation solving, are often the most accessible option even for finding equilibria if the time stepper is already available \cite{dijkstra2014numerical,dijkstra2023bifurcation}. The review \cite{dijkstra2014numerical} refers to several demonstrations of this approach, including bifurcation analysis of stability boundaries and unstable states of a full 3D model of the Atlantic Meridional Overturning Circulation \cite{thies2009bifurcation}. As methods based on time steppers using matrix-free linear algebra only need the existence of an arbitrarily initializable time-stepping algorithm for the dynamical system, they have also been applied as entirely ``\emph{equation-free}'' computations in many-particle simulations in chemistry \cite{kevrekidis2009equation,kevrekidis2010equation} and in agent-based models \cite{thomas2016equation}.

A special case for continuation in PDEs is the tracking of patterns. The general-purpose package \textsc{pde2path} is discussed in Chapter \ref{GohRademacherLloyd24} of this handbook \cite{GohRademacherLloyd24}, making bifurcation analysis possible for general classes of PDEs such as elliptic PDEs or reaction-diffusion equations in $1$ to $3$ space dimensions. In particular, steady or steadily moving PDE patterns (such as traveling or spiral waves) arise as relative equilibria in domains and equations with continuous symmetry. These patterns require a treatment that generalizes the approach in Section~\ref{sec:brus} for discrete symmetries to continuous symmetries.



    
\section{Availability of Data and Access} Scripts reproducing the computational data for all figures can be accessed at \url{https://github.com/jansieber/DS-handbook24-odebifurcations-resources}.

For the purpose of open access, the corresponding author has applied a ‘Creative Commons Attribution (CC BY) licence to any Author Accepted Manuscript version arising from this submission.

\bibliographystyle{ieeetr}
\bibliography{literature_dankowicz}
\end{document}